\documentclass[MSNbibl,number,citesort,secthm,seceqn,dvips]{arxbj}
\usepackage{graphicx}

\aid{0}
\volume{18}
\issue{3}
\pubyear{2012}
\firstpage{1061}
\lastpage{1098}
\doi{10.3150/11-BEJ362}

\makeatletter

\define@key{bid}{pmid}{}

\def\convp{\stackrel{\mathrm{P}}{\rightarrow}}
\def\ex{{\mathbb E}}
\def\var{\operatorname{\mathbb{V}ar}}
\newcommand{\eqref}[1]{(\ref{#1})}

\newtheorem{prop}{Proposition}[section]
\newtheorem{cor}{Corollary}[section]
\newproclaim{cnd}{Condition}[section]
\newtheorem{lemma}{Lemma}[section]
\makeatother

\begin{document}
\begin{frontmatter}

\title{$\sqrt{n}$-consistent parameter estimation for systems of ordinary differential equations: bypassing numerical integration via smoothing}
\runtitle{Parameter estimation for ODEs}

\begin{aug}
\author[1]{\fnms{Shota} \snm{Gugushvili}\corref{}\thanksref{1}\ead[label=e1]{s.gugushvili@vu.nl}} \and
\author[2]{\fnms{Chris A.J.} \snm{Klaassen}\thanksref{2}\ead[label=e2]{c.a.j.klaassen@uva.nl}}
\runauthor{S. Gugushvili and C.A.J. Klaassen}
\address[1]{Department of Mathematics, Vrije Universiteit Amsterdam,
De Boelelaan 1081, 1081 HV Amsterdam, The~Netherlands. \printead{e1}}
\address[2]{Korteweg-de Vries Institute for Mathematics,
Universiteit van Amsterdam, P.O. Box 94248, 1090 GE Amsterdam,
The Netherlands. \printead{e2}}
\end{aug}

\received{\smonth{7} \syear{2010}}
\revised{\smonth{2} \syear{2011}}

%
\begin{abstract}
We consider the problem of parameter estimation for a system of
ordinary differential equations from noisy observations on a
solution of the system. In case the system is nonlinear, as it
typically is in practical applications, an analytic solution to it
usually does not exist. Consequently, straightforward estimation
methods like the ordinary least squares method depend on
repetitive use of numerical integration in order to determine the
solution of the system for each of the parameter values
considered, and to find subsequently the parameter estimate that
minimises the objective function. This induces a huge
computational load to such estimation methods. We study the consistency
of an alternative
estimator that is defined as a minimiser of an appropriate
distance between a nonparametrically estimated derivative of the
solution and the right-hand side of the system applied to a
nonparametrically estimated solution. This smooth and match estimator
(SME) bypasses
numerical integration altogether and reduces the amount of
computational time drastically compared to ordinary least squares.
Moreover, we show that under suitable regularity conditions this
smooth and match estimation procedure leads to a~$\sqrt{n}$-consistent
estimator of
the parameter of interest.
\end{abstract}

%

\end{frontmatter}

\section{Brief introduction}\label{briefintroduction}\vspace*{-3pt}

Many dynamical systems in science and applications are modeled by
a $d$-dimensional system of ordinary differential equations, denoted as
%
\begin{equation}\label{model}
\cases{
{ x}^{\prime}(t)=F({ x}(t),{ \theta}),& \quad $t\in[0,1]$,\cr
x(0)= \xi,
}\vadjust{\goodbreak}
\end{equation}
where $\theta$ is the unknown parameter of interest and $\xi$ is the
initial condition. With
$x_\theta(t)$ the solution vector corresponding to the parameter
value $\theta,$ we observe
\[
Y_{ij}=x_{\theta j}(t_i)+\epsilon_{ij},\qquad
i=1,\ldots,n, j=1,\ldots,d,
\]
where the observation times $0\leq t_1<\cdots <t_n\leq1$ are known
and the random variables~$\epsilon_{ij}$ have mean 0 and model
measurement errors combined with latent random deviations from the
idealised model \eqref{model}. Under regularity conditions, the
ordinary least squares estimator
%
\begin{equation}\label{OLSE}
\tilde{\theta}_n=\operatorname{\arg\min}\limits_\eta\sum_{i=1}^n
\sum_{j=1}^d  \bigl(Y_{ij}-x_{\eta j}(t_i) \bigr)^2
\end{equation}
of $\theta$ is $\sqrt n$-consistent, at least theoretically. For
systems \eqref{model} that do not have explicit solutions, one
typically uses iterative procedures to approximate this ordinary least squares
estimator. However, since every iteration in such a procedure
involves numerical integration of the system \eqref{model} and
since the number of iterations is typically very large, in
practice it is often extremely difficult if not impossible to
compute \eqref{OLSE}, cf. page 172 in \cite{voit}. Here, we
present a feasible and
computationally much faster method to estimate the parameter
$\theta$. To define the estimator of $\theta$, we first construct kernel
estimators
\[
\hat{x}_j(t)=\sum_{i=1}^n (t_i-t_{i-1}) \frac1b
K \biggl(\frac{t-t_i}{b}  \biggr) Y_{ij}
\]
of $x_{\theta j}$ with $K$ a kernel function and $b=b_n$ a
bandwidth. Now, the estimator $\hat{\theta}_n$ of $\theta$ is
defined as
%
\begin{equation}
\label{2half}
\hat{\theta}_n =\operatorname{\arg\min}\limits_\eta\int_0^1 \|
\hat{x}^{\prime}(t) - F(\hat{x}(t),\eta)  \|^2 w(t)\,\mathrm{d}t,
\end{equation}
where $\|\cdot\|$ denotes the usual Euclidean norm and $w(\cdot)$
is a weight function. Related approaches have been suggested in
computational biology and numerical analysis literature, see, for
example, \cite{bellman,savageau} and \cite{varah}.

The main result of this paper is that this smooth and match estimator
$\hat{\theta}_n$ is $\sqrt
n$-consistent under mild regularity conditions. So, asymptotically the
SME $\hat{\theta}_n$
is comparable to the ordinary least squares estimator in
statistical performance, but it avoids the computationally costly
repeated use of numerical integration of \eqref{model}.

\section{Introduction}\label{introduction}

Let us introduce the contents of this paper in more detail.
Systems of ordinary differential equations play a fundamental role
in many branches of natural sciences, for example, mathematical biology,
see \cite{keshet}, biochemistry, see \cite{voit}, or the theory
of chemical reaction networks in general, see, for instance,
\cite{feinberg} and \cite{sontag}. Such systems usually depend
on parameters, which in practice are often only approximately
known, or are plainly unknown. Knowledge of these parameters is
critical for the study of the dynamical system or process that the
system of ordinary differential equations describes. Since these
parameters usually cannot be measured directly, they have to be
inferred from, as a rule, noisy measurements of various quantities
associated with the process under study. More formally, in this
paper we consider the following setting: let, as in \eqref{model},
%
\begin{equation}
\label{system}
\cases{
{ x}^{\prime}(t)=F({ x}(t),{ \theta}), &\quad $t\in[0,1]$,\cr
x(0)=\xi,
}
\end{equation}
be a system of autonomous differential equations depending on a
vector of real-valued parameters. Here ${
x}(t)=(x_1(t),\ldots,x_d(t))^{\mathrm{T}}$ is a $d$-dimensional state
variable, ${ \theta}=(\theta_1,\ldots,\theta_p)^{\mathrm{T}}$ denotes a
$p$-dimensional parameter, while the column $d$-vector $x(0)=\xi$
defines the initial condition. Whether the latter is known or
unknown, is not relevant in the present context, as long as it
stays fixed. Denote a solution to \eqref{system} corresponding to
parameter value $\theta$ by $x_\theta(t)=(x_{\theta
1}(t),\ldots,x_{\theta d}(t))^{\mathrm{T}}.$ Suppose that at known time
instances $0\leq t_1<\cdots <t_n\leq1$ noisy observations
%
\begin{equation}
\label{onehalf}
Y_{ij}=x_{\theta j}(t_i)+\epsilon_{ij}, \qquad
i=1,\ldots,n, j=1,\ldots,d,
\end{equation}
on the solution $x_\theta$ are available. The random variables
$\epsilon_{ij}$ model measurement errors, but they might also
contain latent random deviations from the idealized model
(\ref{model}). Such random deviations are often seen in real-world
applications. Based on these observations, the goal is to infer
the value of $\theta,$ the parameter of interest.

The standard approach to estimation of $\theta$ is based on the
least squares method (the least squares method is credited to Gau{\ss}
and Legendre, see \cite{stigler}), see, for example, \cite{hemker}
and \cite{stortelder}. The least squares estimator is defined as a
minimiser of the sum of squares, that is,
\[
\tilde{\theta}_n=\operatorname{\arg\min}\limits_{\eta}
R_n(\eta)=\operatorname{\arg\min}\limits_{\eta}\sum_{i=1}^n\sum_{j=1}^d
\bigl(Y_{ij}-x_{\eta j}(t_i)\bigr)^2.
\]
If the measurement errors are Gaussian, then $\tilde{\theta}_n$
coincides with the maximum likelihood estimator and is
asymptotically efficient. Since the differential equations setting
is covered by the general theory of nonlinear least squares,
theoretical results available for the latter apply also in the
differential equations setting and we refer for example, to~\cite{jennrich} and~\cite{wu}, or more generally to
\cite{geer1,geer2}, and \cite{pollard} for a thorough
treatment of the asymptotics of the nonlinear least squares
estimator. The paper that explicitly deals with the ordinary
differential equations setting is \cite{xue}. Despite its appealing
theoretical properties, in
practice the performance of the least squares method can
dramatically degrade if \eqref{system} is a nonlinear
high-dimensional system and if $\theta$ is high-dimensional. In
such a~case, we have to face a nonlinear optimisation problem
(quite often with many local minima) and search for a global
minimum of the least squares criterion function~$R_n$ in a~high-dimensional
parameter space. The search process is most often
done via gradient-based methods, for example, the Levenberg--Marquardt
method, see \cite{levenberg}, or via random search algorithms,
see Section 4.5.2 in \cite{voit} for a literature overview. Since
nonlinear systems in general do not have solutions in closed form,
use of numerical integration within a~gradient-based search method
and serious computational time associated with it seem to be
inevitable. For instance, in a relatively simple example of a
four-dimensional system considered in Appendix 1 of
\cite{almeida}, repetitive numerical
integration of the system takes up to 95\% of the total computational
time required to compute the least squares estimator via a~gradient
based optimisation method. Likewise,
random search algorithms are also very costly computationally and in
general, computational time will typically be a~problem for any
optimisation algorithm that relies on numerical integration of any
relatively realistic nonlinear system of ordinary differential
equations, cf. page 172 in \cite{voit}. One example is furnished by
\cite{kikuchi}, where a
system that consists of five differential equations and contains sixty
parameters and that describes a simple gene regulatory network
from~\cite{hlavacek} is considered. The optimisation algorithm (a
genetic algorithm) was run for seven loops each lasting for about
ten hours on the AIST CBRC Magi Cluster with 1040 CPUs (Pentium
III 933 MHz).\footnote{See \url{http://www.cbrc.jp/magi} for the
cluster specifications.} This amounted to a total of ca. 70\,000
CPU hours. The authors also remarked that the gradient-based
search algorithm would not be feasible in their setting at all. The
problems become aggravated for systems of ordinary differential
equations that exhibit stiff behaviour, that is, systems that contain
both `slow' and `fast' variables and that are
difficult to integrate via explicit numerical integration schemes,
see, for example,~\cite{hairer} for a comprehensive treatment of methods
of solving numerically stiff systems. Even if a system is not
stiff for the true parameter value $\theta,$ during the numerical
optimisation procedure one might pass the vicinity of parameters
for which the system is stiff, which will necessarily slow down
the optimisation process.

The Bayesian approach to estimation of $\theta,$ see, for example,
\cite{gelman} and \cite{girolami}, encounters similar huge
computational problems. In the Bayesian approach, one puts a prior
on the parameter $\theta$ and then obtains the posterior via
Bayes' formula. The posterior contains all the information
required in the Bayesian paradigm and can be used to compute for
example,
point estimates of $\theta$ or Bayesian credible intervals. If
$\theta$ is high-dimensional, the posterior will typically not be
manageable by numerical integration and one will have to resort to
Markov Chain Monte Carlo (MCMC) methods. However, sampling from
the posterior distribution for $\theta$ via MCMC necessitates at
each step numerical integration of the system \eqref{system}, in
case the latter does not have a closed form solution.
Computational time might thus become a problem in this case as
well. Also, since in general the likelihood surface will have a
complex shape with many local optima, ripples, and ridges, see,
for example, \cite{girolami} for an example, serious convergence
problems might arise for MCMC samplers.

Yet another point is that in practice both the least squares
method and the Bayesian approach require good initial guesses of
the parameter values. If these are not available, then both approaches
might have problems with convergence to the true
parameter value within a reasonable amount of time.

Over the years a number of improvements upon the classical methods to
compute the least squares estimator have been proposed in the
literature. In particular, the multiple shooting method of \cite
{bock} and the interior-point or barrier method for large-scale
nonlinear programming as in \cite{biegler} have proved to be quite
successful. These two approaches tend to be much more stable than
classical gradient-based methods, have a better chance to converge even
from poor initial guesses of parameters, and in general require a far
less number of iterations until convergence is achieved. However, they
still require a nontrivial amount of computational power.

A general overview of the typical difficulties in parameter estimation
for systems of ordinary differential equations is given in
\cite{ramsay}, to which we refer for more details. For a~recent
overview of typical approaches to parameter estimation for systems
of ordinary differential equations in biochemistry and associated
challenges see, for example, \cite{chou}.

To evade difficulties associated with the least squares method, or
more precisely with numerical integration that it usually
requires, a two-step method was proposed in \cite{bellman} and \cite
{varah}. In the
first step, the solution $x_{\theta}$ of \eqref{system} is
estimated by considering estimation of the individual components
$x_{\theta1},\ldots,x_{\theta d}$ as nonparametric regression
problems and by using the regression spline method for estimation of
these components. The derivatives of $x_{\theta
1},\ldots,x_{\theta d}$ are also estimated from the data by
differentiating the estimators of $x_{\theta1},\ldots,x_{\theta
d}$ with respect to time $t.$ Thus, no numerical integration of the
system \eqref{system} is needed. In the second step, the obtained
estimate of $x_{\theta}$ and its derivative $x_{\theta}^{\prime}$
are plugged into \eqref{system} and an estimator of $\theta$ is
defined as a minimiser in $\theta$ of an appropriate distance
between the estimated left- and right-hand sides of \eqref{system} as
for example, in \eqref{2half}. Since this estimator of $\theta$
results from a minimisation procedure,
it is an M-estimator, see, for example, the
classical monograph \cite{huber}, or Chapter 7 of \cite{bkrw},
Chapter 5 of \cite{vdvaart}, and Chapter~3.2 of \cite{wellner}
for a more modern exposition of the theory of M-estimators. For an
approach to estimation of $\theta$ related to \cite{bellman} and
\cite{varah} see also
\cite{savageau}, as well as \cite{almeida}, where a practical
implementation based on neural networks is studied. The intuitive
idea behind the use of this two-step estimator is clear: among all
functions defined on $[0,1],$ any reasonably defined distance
between the left- and right-hand side of \eqref{system} is minimal
(namely, it is zero) for the solution $x_{\theta}$ of
\eqref{system} and the true parameter value $\theta.$ For
estimates close enough in an appropriate sense to the solution
$x_{\theta},$ the minimisation procedure will produce a minimiser
close to the true parameter value, provided certain
identifiability and continuity conditions hold. This intuitive
idea was exploited in \cite{brunel}, where a more general setting
than the one in \cite{bellman} and \cite{varah} was considered.
Another paper in the
same spirit as \cite{bellman} and \cite{varah} is \cite{liang}.

This two-step approach will typically lead to considerable savings
in computational time, as unlike the straightforward least squares
estimator, in its first step it just requires finding
nonparametric estimates of $x_{\theta}$ and $x_{\theta}^{\prime},$
for which fast and numerically reliable recipes are available,
whereas the gradient-based least squares method will still rely on
successive numerical integrations of \eqref{system} for different
parameter values $\theta$ in order to find a global minimiser
minimising the least squares criterion function. We refer to
\cite{almeida} for a particular example demonstrating gains in
the computational time achieved by the two-step estimator in
comparison to the ordinary least squares estimator. When the
right-hand side $F$ of \eqref{system} is linear in
$\theta_1,\ldots,\theta_p$ and $d=1,$ further simplifications will
occur in
the second step of the two-step estimation procedure, as one will
essentially only have to face a weighted linear regression problem
then. This is unlike the least squares approach, which cannot
exploit linearity of $F$ in $\theta_1,\ldots,\theta_p.$ However,
we would also like to stress the fact that the two-step estimator
does not necessarily have to be considered a~competitor of either
the least squares or the Bayesian approach. Indeed, since in
practice both of these approaches require good initial guesses for
parameter values, these can be supplied by the two-step estimator.
In this sense, the proposed two-step estimation approach can be
thought of as complementing both the least squares and the
Bayesian approaches. Moreover, an additional modified
Newton--Raphson step suffices to arrive at an estimator that is
asymptotically equivalent to the exact ordinary least squares
estimator, as will be shown elsewhere.

A certain limitation of the two-step approach is that it requires that
measurements on all state variables $x_{\theta j},j=1,\ldots,d$, are
available. The latter is not always the case in practical applications.
In some cases, the unobserved variables can be eliminated by
transforming the first order system into a higher order one and next
applying a generalisation of the smooth and match method to this higher
order system. Under appropriate conditions, this approach should yield
a consistent estimator. Although one can always formally perform the
least squares procedure, without further assumptions on the system it
is far from clear that it leads to a consistent estimator of the
parameter of interest.

Our goal in the present work is to undertake a rigorous study of the
asymptotics of a two-step estimator of $\theta.$ Our exposition is
similar to that in
\cite{brunel} to some degree, but one of the differences is that
instead of regression spline estimators we use kernel-type estimators for
estimation of $x_{\theta}$ and $x_{\theta}^{\prime}.$\footnote{The
proofs of the main results in \cite{brunel} are incomplete and
the main theorems require further conditions in order to
hold.} The conditions are also different. We hope that
our contribution will motivate further research into the
interesting topic of parameter estimation for systems of ordinary
differential equations.

There exists an alternative approach to
the ones described here, which also employs nonparametric smoothing,
see \cite{ramsay}. For information on its asymptotic properties, we
refer to \cite{qi}. For nonlinear systems, this approach will
typically reduce to one of the realisations of the ordinary least
squares method, for example, Newton--Raphson algorithm, where however
numerical integration of \eqref{system} will be replaced by
approximation of the solution of the system \eqref{system} by an
appropriately chosen element of some finite-dimensional function space.
This seems to reduce considerably the computational load in comparison
to the gradient-based optimisation methods which employ numerical
integration of \eqref{system}. However, it still appears to be
computationally more intense than the two-step approach advocated in
the present work.

We conclude the discussion in this section by noting that when modeling
various processes, some authors prefer not to specify the right-hand
side of \eqref{system} explicitly (the latter amounts to explicit
specification of the $F(\cdot,\cdot)$ in \eqref{system}), but simply
assume that the right-hand side of \eqref{system} is some unknown
function of $x,$ that is, is given by $F(x(t))$ with $F$ unknown, and
proceed to its estimation via nonparametric methods, see, for example,~\cite{ellner}.
This has an advantage of safeguarding against possible
model misspecification. However, the question whether one has or has
not to specify $F$ explicitly appears to us to be more of a
philosophical nature and boils down to a discussion on the use of
parametric or nonparametric models, that is, whether one has strong
enough reasons to believe that the process under study can be described
by a model as in \eqref{system} with $F$ known or not. We do not
address this question here, because an answer to it obviously depends
on the process under study and varies from case to case. For a related
discussion, see \cite{hooker}.

The rest of the paper is organised as follows: in the next section,
we will detail the approach that we use and present its
theoretical properties. In particular, we will show that under
appropriate conditions our two-step approach leads to a consistent
estimator with a $\sqrt{n}$ convergence rate, which is the best
possible rate in regular parametric models.\footnote{
It is claimed in \cite{liang} that the two-step estimation procedure suggested
there leads to a faster rate than $\sqrt{n},$ which is impossible.
Indeed, Theorem 2 of \cite{liang} and its proof are incorrect.} Section \ref{discussion}
contains a discussion on the obtained results together with simulation
examples. The proofs of the main results are relegated to
Section \ref{proofs}, while Appendix \ref{apps} contains some auxiliary
statements.

\section{Results}
\label{results}

First of all, we point out that in the present study we will be
concerned with the asymptotic behaviour of an appropriate two-step
estimator of $\theta$ under a suitable sampling scheme. We will
primarily be interested in intuitively understanding the
behaviour of a~relatively simple estimator of $\theta,$ as well as in
a clear presentation of the obtained results and the proofs.
Consequently, the stated conditions will not always be minimal and
can typically be relaxed at appropriate places.

We first define the sampling scheme.

\begin{cnd}
\label{cnd_times} The observation times $0\leq t_1<\cdots <t_n\leq
1$ are deterministic and known and there exists a constant $c_0\geq1,$ such
that for all $n$
\[
\max_{2\leq i\leq n} |t_i-t_{i-1} |\leq\frac{c_0}{n}
\]
holds. Furthermore, there exists a constant $c_1\geq1,$ such that for
any interval $A\subseteq[0,1]$ of length $|A|$ and all $n\geq1$
the inequality
\[
\frac{1}{n}\sum_{i=1}^n 1_{[t_i\in A]}\leq c_1
\max \biggl(|A|,\frac{1}{n} \biggr)
\]
holds.
\end{cnd}

Hence, we observe the solution of the system
\eqref{system} on the interval $[0,1].$ Instead of $[0,1]$ we
could have taken any other bounded interval. Conditions on
$t_1,\ldots,t_n$ as in Condition~\ref{cnd_times} are typical in
nonparametric regression, see, for example, \cite{gasser} and
Section~1.7
in \cite{tsybakov}, and they imply that $t_1,\ldots,t_n$ are
distributed over $[0,1]$ in a sufficiently uniform manner. The
most important example in which Condition \ref{cnd_times} is
satisfied, is when the observations are spaced equidistantly over
$[0,1],$ that is, when $t_j=j/n$ for $j=1,\ldots,d.$ In this case, one
may take $c_0=c_1=2.$ Notice that we do not necessarily assume that the initial
condition $x(0)=\xi$ is measured or is known. If it is, then it is
incorporated into the observations and is used in the first step
of the two-step estimation procedure.

\begin{cnd}
\label{cnd_epsilon} The random variables $\epsilon_{ij} ,
i=1,\ldots,n, j=1,\ldots,d,$ from \eqref{onehalf} are independent
and are normally distributed with mean zero and finite variance
$\sigma_j^2.$
\end{cnd}

This assumption of Gaussianity of the $\epsilon_{ij}$'s
may be dropped in various ways, as we will see below; see the note
after Proposition \ref{priestleyprop} and Appendix \ref{app2}.

We next state a condition on the parameter set.

\begin{cnd}
\label{cnd_theta} The parameter set $\Theta$ is a compact subset
of ${\mathbb{R}}^p.$
\end{cnd}

Compactness of $\Theta$ allows one to put relatively weak
conditions on the structure of the system \eqref{system}, that is,
the function $F.$

Just as the least squares method, see, for example, \cite{jennrich},
our smooth and match approach also requires some regularity of the solutions
of \eqref{system}. In what follows, a derivative of any function
$f$ with respect to the variable $y$ will be denoted by
$f_y^{\prime}.$ For the second derivative of $f$ with respect to
$y$, we will use the notation $f_{yy}^{\prime\prime}$ with a
similar convention for mixed derivatives. An integral of a vector- or
matrix-valued function will be understood componentwise.

\begin{cnd}
\label{cnd_solution} The following conditions hold:
\begin{longlist}[(iii)]
\item[(i)]
the mapping $F\dvtx  {\mathbb{R}}^d \times\Theta\rightarrow
{\mathbb{R}}^d$ from \eqref{system} is such that its second
derivatives~$F_{\theta\theta}^{\prime\prime},F_{\theta
x}^{\prime\prime},F_{xx}^{\prime\prime}$ are continuous;
\item[(ii)]
for all parameter values $\theta\in\Theta,$ the
solution $x_{\theta}$ of \eqref{system} is defined on the interval $[0,1];$
\item[(iii)]
for all parameter values $\theta\in\Theta,$ the
solution $x_{\theta}$ of \eqref{system} is unique on $[0,1];$
\item[(iv)]
for all parameter values $\theta\in\Theta,$ the
solution $x_{\theta}$ of \eqref{system} is a $C^{\alpha}$ function
of $t$ on the interval $[0,1]$ for some
positive integer $\alpha.$
\end{longlist}
\end{cnd}

Observe that Condition \ref{cnd_solution}(i) implies
existence and uniqueness of the solution of \eqref{system} in some
neighbourhood of $0.$ However, we want the existence and
uniqueness to hold on the whole interval $[0,1]$ and therefore a
priori require (ii) and (iii). Furthermore, $\alpha\geq2$ in (iv)
is required when establishing appropriate asymptotic properties of
nonparametric estimators of the solution $x_{\theta}$ and its
derivative, while $\alpha\geq3$ is needed in Propositions~\ref{functional}
and \ref{linfunct}, and $\alpha\geq4$ in
Theorem~\ref{mainthm}, respectively. Notice that for every
$\theta$ the solution~$x_{\theta}$ is of class $C^{\alpha}$ in $t$
in a neighbourhood of $0,$ provided for a given $\theta$ the
function $F$ is of class $C^{\alpha}$ in its first argument.
However, we want this to hold on the whole interval $[0,1]$ and
therefore require (iv). Since in the theory of chemical reaction
networks, see, for instance, \cite{sontag}, the components of $F$
are usually polynomial or rational functions of
$x_1,\ldots,x_d$
and $\theta_1,\ldots,\theta_p,$ the solution of \eqref{system}
will be smooth enough in many examples and $\alpha\geq4$ is
satisfied in a large number of practical examples. For the
above-mentioned facts from the theory of ordinary differential
equations, see, for example, Chapter 2 in \cite{arnold}. Also notice that
the condition on $F$ in \cite{liang}, see Assumption C on page
1573, puts severe restrictions on $F$ and excludes for example, quadratic
nonlinearities of $F$ in $x_1,\ldots,x_d.$ This, of course, has to
be avoided.

Recall that our observations are $Y_{ij}=x_{\theta
j}(t_i)+\epsilon_{ij}$ for $i=1,\ldots,n,j=1,\ldots,d.$ We propose
the following nonparametric estimator for $x_{\theta j},$
%
\begin{equation}
\label{estx} \hat{x}_j(t)=\sum_{i=1}^n (t_i-t_{i-1})
\frac{1}{b}K \biggl(\frac{t-t_i}{b}  \biggr) Y_{ij},
\end{equation}
where $K$ is a kernel function, while the number $b=b_n>0$ denotes
a bandwidth that we take to depend on the sample size $n$ in such
a way that $b_n\rightarrow0$ as $n\rightarrow\infty.$ In line
with a traditional convention in kernel estimation theory, we will
suppress the dependence of~$b_n$ on $n$ in our notation, since no
confusion will arise. When the $t_i$'s are equispaced, the
estimator \eqref{estx} can in essence be obtained by modifying the
Nadaraya--Watson regression estimator, cf. page 34 in
\cite{tsybakov}. It is usually called the Priestley--Chao
estimator after the authors who first proposed it in
\cite{priestley}. As far as an estimator of $x_{\theta
j}^{\prime}(t)$ is concerned, we define it as the derivative of
$\hat{x}_j(t)$ with respect to $t,$ choosing $K$ as a
differentiable function. Notice that the bandwidth $b$ plays a
role of regularisation parameter: too small a bandwidth results in
an estimator with small bias, but large variance, while too large
a~bandwidth results in an estimator with small variance, but large
bias, see, for example, pages 7--8 and 32 in \cite{tsybakov} for a
relevant discussion. In principle one could use different
bandwidth sequences for estimation of $x_j$ for different $j$'s,
but as can be seen from the proofs in Section \ref{proofs},
asymptotically this will not make a difference for an estimator of
$\theta$. A similar remark applies to the use of different
bandwidths for estimation of $x_{\theta j}$ and its derivative
$x_{\theta j}^{\prime}.$ Arguably, the estimator \eqref{estx} is
simple and there exist other estimators that may outperform it in
certain respects in practice. However, as we will show later on,
even such a simple estimator leads to a $\sqrt{n}$-consistent
estimator of $\theta.$

Theoretical properties of the Priestley--Chao estimator were
studied in \cite{benedetti,priestley,schuster}.
However, the first two papers do not cover its convergence in the
$L_{\infty}$ (supremum) norm, while the third one does not do it
in the form required in the present work. Since this is needed in
the sequel, we will supply the required statement, see Proposition
\ref{priestleyprop} below.

To put things in a somewhat more
general context than the one in our differential equations
setting, consider the following regression model:
%
\begin{eqnarray}
\label{regression}
&&Y_i=\mu(t_i)+\epsilon_i,\qquad i=1,\ldots,n,\nonumber\\[-2pt]
&&t_1,\ldots,t_n \mbox{ satisfy Condition \ref{cnd_times}},\\[-2pt]
&&\epsilon_1,\ldots,\epsilon_n \mbox{ are i.i.d. Gaussian with }
\ex[\epsilon_i]=0 \mbox{ and } \var[\epsilon_i]=\sigma^2>0.\nonumber
\end{eqnarray}
Our goal is to estimate the regression function $\mu$ and its
derivative~$\mu^{\prime}.$ The estimator of~$\mu$ will be given by
an expression similar to \eqref{estx}, namely
%
\begin{equation}
\label{priestley} \hat{\mu}_n(t)=\sum_{i=1}^n (t_i-t_{i-1}) \frac
1b K \biggl(\frac{t-t_i}{b}  \biggr) Y_i ,\vadjust{\goodbreak}
\end{equation}
while an estimator of $\mu^{\prime}$ will be given by
$\hat{\mu}_n^{\prime}.$ We postulate the following condition on
the kernel $K$ for some strictly positive integer $\alpha.$
\begin{cnd}
\label{cndkernel} The kernel $K$ is symmetric and twice
continuously differentiable, it has support within $[-1,1],$ and
it satisfies the integrability conditions: $\int_{-1}^1 K(u)\,\mathrm{d}u=1$
and $\int_{-1}^1 u^{\ell}K(u)\,\mathrm{d}u=0$ for $\ell=1,\ldots,\alpha-1.$
If $\alpha=1,$ only the first of the two integrability conditions is required.
\end{cnd}

The following proposition holds.

\begin{prop}
\label{priestleyprop} Suppose the regression model
\eqref{regression} is given and Condition \ref{cndkernel} holds.
Fix a number $\delta,$ such that $0<\delta<1/2.$
\begin{longlist}[(ii)]
\item[(i)] If $\mu$ is $\alpha\geq1$ times continuously differentiable
and $b\rightarrow0$ as $n\rightarrow\infty,$ then
%
\begin{equation}
\label{mup} \sup_{t\in[\delta,1-\delta]}|\hat{\mu}_n(t)-\mu(t)| =
\mathrm{O}_P \Biggl(
b^{\alpha}+\frac{1}{nb^2}+\sqrt{\frac{\log n}{nb}}
 \Biggr).
\end{equation}

\item[(ii)] If $\mu$ is $\alpha\geq2$ times continuously differentiable and
$b\rightarrow0$ as $n\rightarrow\infty,$
then
%
\begin{equation}
\label{muprimep}
\sup_{t\in[\delta,1-\delta]}|\hat{\mu}_n^{\prime}(t)-\mu
^{\prime}(t)|
= \mathrm{O}_P \Biggl(
b^{\alpha-1}+\frac{1}{nb^3}+\sqrt{\frac{\log
n}{nb^3}} \Biggr)
\end{equation}
is valid. In particular, $\hat{\mu}_n$ and $\hat{\mu}_n^{\prime}$
are consistent on $[\delta,1-\delta],$ if $nb^3/\log
n\rightarrow\infty$ holds additionally.
\end{longlist}
\end{prop}

Gaussianity of the $\epsilon_i$'s allows one to prove \eqref{mup}
and \eqref{muprimep} by relatively elementary means. This
assumption can be modified in various ways, for instance by
assuming that the $\epsilon_i$'s are bounded, and we state and
prove the corresponding modification of Proposition
\ref{priestleyprop} in Appendix \ref{app2}, see Proposition
\ref{priestleypropbounded}. In general, normality of the
measurement errors is a standard assumption in parameter
estimation for systems of ordinary differential equations, see,
for example, \cite{girolami,hemker}, and \cite{ramsay}.

The following corollary is immediate from Proposition
\ref{priestleyprop}.
\begin{cor}
\label{corxi} Let $\alpha$ be the same as in Condition \ref
{cnd_solution}. Under Conditions \ref{cnd_times}--\ref{cndkernel}, we
have for the estimator $\hat{x}_j$
%
\begin{equation}
\label{xip} \sup_{t\in[\delta,1-\delta]}|\hat{x}_j(t)-x_{\theta
j}(t)| = \mathrm{O}_P\Biggl(
b^{\alpha}+\frac{1}{nb^2}+\sqrt{\frac{\log n}{nb}}  \Biggr)
\end{equation}
and
%
\begin{equation}
\label{xiprimep}
\sup_{t\in[\delta,1-\delta]}|\hat{x}_j^{\prime}(t)-x_{\theta
j}^{\prime}(t)| = \mathrm{O}_P\Biggl(
b^{\alpha-1}+\frac{1}{nb^3}+\sqrt{\frac{\log
n}{nb^3}}  \Biggr),
\end{equation}
provided $\alpha\geq2$ and $b\rightarrow0$ as
$n\rightarrow\infty.$ In particular, $\hat{x}_j$ and
$\hat{x}_j^{\prime}$ are consistent, if $nb^3/\log
n\rightarrow\infty$ holds additionally.
\end{cor}

In the proof of Proposition \ref{consistency}, we will apply the
continuous mapping theorem in order to prove convergence in
probability of certain integrals of $F$ and its derivatives with
$\hat{x}_j$'s plugged in. This is where Corollary \ref{corxi} is used.


Now that we have consistent (in an appropriate sense) estimators
of $x_{\theta j}$ and $x_{\theta j}^{\prime},$ from the smoothing step
we can move to the
matching step in the construction of our smooth and match estimator of
$\theta.$ In particular, we define the estimator $\hat{\theta}_n$
of $\theta$ as
%
\begin{eqnarray}
\label{thetan}
\hat{\theta}_n & =&\operatorname{\arg\min}\limits_{\eta\in\Theta} \int_0^1
\|
\hat{x}^{\prime}(t) - F(\hat{x}(t),\eta) \|^2 w(t)\,\mathrm{d}t\nonumber
\\[-8pt]
\\[-8pt]
& =& \operatorname{\arg\min}\limits_{\eta\in\Theta} M_{n,w}(\eta),
\nonumber
\end{eqnarray}
where $\|\cdot\|$ denotes the usual Euclidean norm and $w$ is a
weight function. We will refer to $M_{n,w}(\eta)$ as a (random)
criterion function. Since $\Theta$ is compact and $M_{n,w}$ under our
conditions is
continuous in $\eta,$ the minimiser $\hat{\theta}_n$ always
exists. The fact that $\hat{\theta}_n$ is a~measurable function of
the observations $Y_{ij}$ follows from Lemma 2 of
\cite{jennrich}. Notice that in~\cite{liang} and \cite{varah}
the criterion function is given by
\[
\sum_{i=1}^n \| \widetilde{x}^{\prime}(t_i) -
F(\widetilde{x}(t_i),\eta)\|^2,
\]
where $\widetilde{x}$ and $\widetilde{x}^{\prime}$ are appropriate
estimators of $x_{\theta}$ and $x_{\theta}^{\prime}.$ However, in
order to obtain a $\sqrt{n}$-consistent estimator of $\theta,$ it
is important to use an integral type criterion: the nonparametric
estimators of $x_{\theta}$ and $x_{\theta}^{\prime}$ have a slower
convergence rate than $\sqrt{n}$ and this is counterbalanced by the
integral criterion from \eqref{thetan}. Indeed, stationarity at $\hat
{\theta}_n$ leads to \eqref{*}. The first factor at the left-hand
side of this equality converges to a constant nondegenerate matrix and
the right-hand side behaves like a linear combination of the
observations with coefficients of order $1/n$ thanks to the
integration; cf. Proposition \ref{linfunct} and its proof. In
light of this the choice of the weight function $w$ also appears
to be important. Furthermore, the observations $Y_{ij}$ from \eqref
{onehalf} indirectly
carry information on the entire curves $x_{\theta
j}(t),t\in[0,1],$ and not only on the points $x_{\theta j}(t_i)$.
An integral type criterion allows one to exploit this fact in the
second step of this smooth and match procedure.

Introduce the asymptotic criterion
\[
M_w(\eta)= \int_0^1 \| F(x_{\theta}(t),\theta) -
F(x_{\theta}(t),\eta) \|^2w(t)\,\mathrm{d}t
\]
corresponding to $M_{n,w}.$ Observe that by Condition
\ref{cnd_solution} it is bounded. Using Corollary~\ref{corxi} as a
building block, one can show that the SME
$\hat{\theta}_n$ is consistent. To this end, we will need the
following condition on the weight function $w$.

\begin{cnd}
\label{cndw} The weight function $w$ is a nonnegative function
that is continuously differentiable, is supported on the interval
$(\delta,1-\delta)$ for some fixed number $\delta,$ such that
$0<\delta<1/2,$ and is such that the
Lebesgue measure of the set $\{t\dvt w(t)>0\}$ is positive.
\end{cnd}

The fact that $w$ vanishes at the endpoints of the interval
$[\delta,1-\delta]$ and beyond, is needed to obtain a
$\sqrt{n}$-consistent estimator of $\theta.$ In particular, together
with differentiability of $w$ it is used in order to establish \eqref
{6*}. The
condition that $w$ is supported on $(\delta,1-\delta)$ takes care
of the boundary bias effects characteristic of the conventional
kernel-type estimators, see, for example, \cite{gasser} for more
information on this. Boundary effects in kernel estimation are
usually remedied by using special boundary kernels, see, for example,
\cite{vanes,gasser2,goldstein2}. Using such a
kernel, it can be expected that in our case as well the boundary
effects will be eliminated and one may relax the requirement
$0<\delta<1/2$ from Condition \ref{cndw} to $\delta=0,$
that is, to allowing $w$ to be supported on $(0,1).$ The condition
that the weight function $w$ is positive on a set with positive
Lebesgue measure, is important for \eqref{identifiability} to hold
and in fact $w(t)=0$ a.e. would be a strange choice.

The following proposition is valid.

\begin{prop}
\label{consistency} Suppose $b\rightarrow0$ and $nb^3/\log n
\rightarrow\infty.$ Under Conditions \ref{cnd_times}--\ref{cndw}
and the additional identifiability condition
%
\begin{equation}
\label{identifiability} \forall\varepsilon>0 \qquad
\inf_{\|\eta-\theta\|\geq\varepsilon} M_w(\eta)>M_w(\theta),
\end{equation}
we have $\hat{\theta}_n\convp\theta.$
\end{prop}

The proposition is proved via a reasoning standard in the theory
of M-estimation: we show that $M_{n,w}$ converges to $M_w$ and
that the convergence is strong enough to imply the convergence of
a minimiser $\hat{\theta}_n$ of $M_{n,w}$ to a minimiser $\theta$
of $M_w,$ cf. Section~5.2 of \cite{vdvaart}. A necessary
condition for \eqref{identifiability} to hold is that
$x_{\theta}(\cdot)\neq x_{\theta^{\prime}}(\cdot)$ for $\theta
\neq
\theta^{\prime}.$ The latter is a minimal assumption for
the statistical identifiability of the parameter $\theta.$ The
identifiability condition \eqref{identifiability} is common in the
theory of M-estimation, see Theorem~5.7 of~\cite{vdvaart}. It
means that $\theta$ is a point of minimum of $M_{w}(\eta)$ and
that it is a {\emph{well-separated}} point of minimum. The most
trivial example with this condition satisfied is when $d=p=1$ and
$x^{\prime}(t)=\theta x(t)$ hold with initial condition
$x(0)=\xi,$ where $\xi\neq0.$ In fact, in this case
\[
M_w(\eta)=(\theta-\eta)^2 \xi^2 \int_{\delta}^{1-\delta}
\mathrm{e}^{2\theta t} w(t)\,\mathrm{d}t,
\]
and this is zero for $\eta=\theta$ and is strictly positive for $\eta
\neq\theta,$ whence \eqref{identifiability} follows. More generally,
since $\Theta$ is
compact and $M_w$ is continuous, uniqueness of a minimiser of
$M_w$ will imply \eqref{identifiability}, cf. Exercise 27 on page 84
of \cite{vdvaart}.

In practice, \eqref{identifiability} might be difficult to check
globally and one might prefer to concentrate on a simpler local
condition: if the first order condition
$[\mathrm{d}M_w(\eta)/\mathrm{d}\eta]_{\eta=\theta}=0$ holds and if the Hessian
matrix $H(\eta)=(\partial^2 M_w(\eta)/\partial\eta_i\,\partial
\eta_j )_{i,j}$ of $M_w$ is strictly positive definite at
$\theta,$ then \eqref{identifiability} will be satisfied for
$\eta\in\Theta$ restricted to some neighbourhood of $\theta,$
because $M_w$ will have a local minimum at such $\theta$ and a
neighbourhood around it can be taken to be compact with small
enough diameter, so that \eqref{identifiability} holds for $\eta$
restricted to this neighbourhood. The conclusion of the theorem
will then hold for the parameter set restricted to this
neighbourhood of $\theta.$

In a statement analogous to Proposition \ref{consistency},
\cite{brunel} requires that the solutions of \eqref{system}
belong to a compact set ${\mathcal{K}}$ for all $\theta$ and $t$
and that $F$ from (\ref{model}) is Lipschitz in its first argument
$x$ for $x$ restricted to this compact ${\mathcal{K}}$ uniformly
in $\theta\in\Theta.$ It is also assumed that the nonparametric
estimators $\hat{x}_n(t)$ belong a.s. to ${\mathcal{K}}$ for all
$n$ and $t.$ However, the latter typically will not hold for
linear smoothers, see Definition 1.7 in \cite{tsybakov}, which
constitute the most popular choice of nonparametric regression
estimators in practice. For instance, local polynomial estimators,
see Section 1.6 in \cite{tsybakov}, projection estimators, see
Section 1.7 in \cite{tsybakov}, or the Gasser--M\"uller estimator,
see \cite{gasser}, are all examples of linear smoothers. Hence, we
prefer to avoid this condition altogether, although this somewhat
complicates the proof.

Under the conditions in this section, it turns out that the
estimator $\hat{\theta}_n$ is not merely a consistent estimator,
but a $\sqrt{n}$-consistent estimator of $\theta,$ in the sense of
\eqref{rootn} below. This result follows in essence from the fact
that up to a higher order term the difference $\hat{\theta}_n
-\theta$ can be represented as the difference of the images of
$\hat{x}$ and $x_{\theta}$ under a certain linear mapping, cf.
Proposition \ref{functional}. It is known that even though
nonparametric curve
estimators cannot usually attain the $\sqrt{n}$ convergence rate,
see, for example, Chapters~1 and 2 of \cite{tsybakov}, extra smoothness
often coming from the structure of linear functionals allows one
to construct in many cases $\sqrt{n}$-consistent estimators of
these functionals via plugging in nonparametric estimators, see,
for example, \cite{bickel} and \cite{goldstein} for more information.
The variance of such plug-in estimators can often be proven to be
of order $n^{-1},$ while the squared bias can be made of order
$n^{-1}$ by undersmoothing, that is, selecting the smoothing
parameter smaller than what is an optimal choice in nonparametric
curve estimation when the object of interest is a curve itself,
cf. \cite{goldstein}. Precisely, this happens in our case as
well: if the mean integrated squared error is used as a
performance criterion of a nonparametric estimator, then under our
conditions the optimal bandwidth for estimation of $x_{\theta}$ is
of order $n^{-1/(2\alpha+1)},$ whereas the optimal bandwidth for
estimation of $\theta$ is in fact smaller, see Theorem
\ref{mainthm} below. Note that undersmoothing is a different approach than
the one in \cite{bickel}, where it is assumed that nonparametric
estimators attain the minimax rate of convergence and the
$\sqrt{n}$-rate for estimation of a functional in concrete
examples, if possible, is achieved by different means exploiting extra
smoothness coming from the structure of a functional, see, for
example,
the first example in Section 2 there. In many cases, it can be
proved that such plug-in type estimators are efficient, see
\cite{bickel}. Notice, however, that in our case this will not
imply that $\hat{\theta}_n$ is efficient.

First, we will provide an asymptotic representation for
the difference $\hat{\theta}_n -\theta.$

\begin{prop}
\label{functional} Let $\theta$ be an interior point of $\Theta.$
Suppose that the conditions of Proposition \ref{consistency} hold
and let the matrix $J_{\theta}$ defined by
%
\begin{equation}
\label{brunel5} J_{\theta}=\int_{\delta}^{1-\delta}
(F_{\theta}^{\prime}(x_{\theta}(t),\theta))^{\mathrm{T}}
F_{\theta}^{\prime}(x_{\theta}(t),\theta) w(t)\,\mathrm{d}t
\end{equation}
be nonsingular. Fix $\alpha\geq3.$ If $b \asymp n^{-\gamma}$
holds for $1/(4\alpha-4)<\gamma<1/6,$ then
%
\begin{equation}
\label{brunel3}
\hat{\theta}_n-\theta=\mathrm{O}_P\bigl(J^{-1}_{\theta}\bigl(\Gamma (\hat
{x} )-\Gamma(x_{\theta})\bigr) \bigr)+\mathrm{o}_P(n^{-1/2})
\end{equation}
is valid with the mapping $\Gamma$ given by
%
\begin{equation}
\label{brunel4} \Gamma(z)=\int_{\delta}^{1-\delta}
 \biggl\{-(F_{\theta}^{\prime}(x_{\theta}(t),\theta))^{\mathrm{T}}
F_{x}^{\prime}(x_{\theta}(t),\theta)w(t) -
\frac{\mathrm{d}}{\mathrm{d}t}[(F_{\theta}^{\prime}(x_{\theta}(t),\theta))^{\mathrm{T}}
w(t)] \biggr\}z(t)\,\mathrm{d}t.
\end{equation}
\end{prop}

With the above result in mind, in order to complete the study of
the asymptotics of $\hat{\theta}_n,$ it remains to study the
mapping $\Gamma.$ Clearly, it suffices to study the asymptotic
behaviour of
\[
\Delta(\hat{\mu}_n)-\Delta(\mu)=\int_{\mathbb{R}} v(t)k(t)\hat
{\mu}_n(t)\,\mathrm{d}t-\int_{\mathbb{R}} v(t)k(t)\mu(t)\,\mathrm{d}t,
\]
where $v$ is a known function that satisfies appropriate
assumptions, while $k$ stands either for $w$ or its derivative
$w^{\prime}.$ The next proposition deals with the asymptotics of
$\Delta(\hat{\mu}_n)-\Delta(\mu).$

\begin{prop}
\label{linfunct} Under Conditions \ref{cndkernel} and \ref{cndw}
and for any continuous function $v$, it holds in the regression
model \eqref{regression} that
\[
\Delta(\hat{\mu}_n)-\Delta(\mu)=\mathrm{O}_P(n^{-1/2}),
\]
provided $\mu$ is $\alpha\geq3$ times differentiable and the
bandwidth $b$ is chosen such that $b\asymp n^{-\gamma}$ holds for
$1/(2\alpha)\leq\gamma\leq1/4.$
\end{prop}

Our main result is a simple consequence of Propositions
\ref{functional} and \ref{linfunct}.

\begin{thm}
\label{mainthm} Let $\theta$ be an interior point of $\Theta.$
Assume that Conditions \ref{cnd_times}--\ref{cndw} together with
\eqref{identifiability} hold and that \eqref{brunel5} is
nonsingular. Fix $\alpha\geq4.$ If the bandwidth $b$ is such that
$b\asymp n^{-\gamma}$ holds for $1/(2\alpha)<\gamma<1/6,$ then
%
\begin{equation}
\label{rootn}
\sqrt{n}(\hat{\theta}_n-\theta)=\mathrm{O}_P(1)
\end{equation}
is valid.
\end{thm}

Thus any bandwidth sequences satisfying the conditions in Theorem
\ref{mainthm} are optimal, in the sense that they lead to
estimators of $\theta$ with similar asymptotic behaviour. In
particular, each
of such bandwidth sequences ensures a $\sqrt{n}$ convergence rate
of $\hat{\theta}_n.$ Consequently, dependence of the asymptotic
properties of the estimator $\hat{\theta}_n$ on the bandwidth is
less critical than it typically is in nonparametric curve
estimation. Notice that the condition $\alpha\geq4$ in Theorem
\ref{mainthm} is needed in order to make the conditions in
Propositions~\ref{functional} and~\ref{linfunct} compatible.

\section{Discussion}
\label{discussion}

The main result of the paper, Theorem \ref{mainthm}, is that under
certain conditions for systems of ordinary differential equations
parameter estimation at the $\sqrt{n}$ rate is possible
\emph{without} employing numerical integration. Although we have
shown this in the case when in the first step of the two-step
procedure a particular kernel-type estimator is used, it may be
expected that a similar result holds for other nonparametric
estimators. For instance, the arguments for the Nadaraya--Watson
estimator seem to be similar, with extra technicalities arising for
example, from the fact that it is a ratio of two functions.
Furthermore, from formula \eqref{6*} it can be seen that the proof of
Proposition \ref{functional} requires that the derivative of an
estimator of $x_{\theta}$ be used as an estimator of $x_{\theta
}^{\prime}.$ Not all popular nonparametric estimators of the
derivatives of a regression function are of this type. In practice for
small or moderate sample sizes it
might be advantageous to use more sophisticated nonparametric
estimators than the Priestley--Chao estimator, but asymptotically
this does not make a difference.

Once a $\sqrt{n}$-consistent estimator $\hat{\theta}_n$ of
$\theta$ is available, one might ask for more, namely if one can
construct an estimator that is asymptotically equivalent to the
ordinary least squares estimator \eqref{OLSE} or that is
semiparametrically efficient. It is expected that this can be
achieved without repeated numerical integration of \eqref{model}
by using $\hat{\theta}_n$ as a starting point and performing a
one-step Newton--Raphson type procedure; see, for example, Section 7.8 of
\cite{bkrw} or Chapter 25 of \cite{vdvaart}. We intend to
address this issue of efficient and ordinary least squares
estimation in a separate publication.

Doubtless, the main challenge in implementing the smooth and match
estimation procedure lies in selecting the smoothing parameter
$b.$ This is true for any two-step parameter estimation procedure for
ordinary differential equations, for example, the one based
on the regression splines as in \cite{brunel} or the local
polynomial estimator as in \cite{liang}, and not only for our
specific estimator. Observations that we supply below apply in
principle to any two-step estimator and not only to the specific
kernel-type one considered in the present work. Hence, they are of general
interest.

Some attention has been paid in the literature to the selection of
the smoothing parameter in the context of parameter estimation for
ordinary differential equations. The considered options range from
subjective choices and smoothing by hand to more advanced
possibilities. Perhaps the simplest solution would be to assume
that the targets of the estimation procedure are $x_{\theta j},
j=1,\ldots,d,$ and to select $b$ (a different one for every
component $x_{\theta j}$) via a cross-validation procedure, see,
for example, Section 5.3 in \cite{wasserman} for a description of
cross-validation techniques in the context of nonparametric
regression. This should produce reasonable results, at least for
relatively large sample sizes, cf. simulation examples considered
in \cite{brunel}. However, it is clear from Theorem \ref{mainthm}
and its proof
that despite its simplicity, such a choice of $b$ will be
suboptimal. Another practical approach to bandwidth selection is
computation of $\hat{\theta}_n=\hat{\theta}_n(b)$ for a range of
values of the bandwidth $b$ on some discrete grid $B$ and then choosing
\[
\hat{b}=\operatorname{\arg\min}\limits_{b\in B} \sum_{i=1}^n\sum_{j=1}^d
\bigl(Y_{ij}-x_{\hat{\theta}_n(b)j}(t_i)\bigr)^2.
\]
This seems a reasonable choice, although the asymptotics of $\hat
{\theta}_n(\hat{b})$ are unclear. One other possibility for practical
bandwidth
selection is nothing else but a variation on the plug-in bandwidth
selection method as described for example, in \cite{jones}: one can see
from the proof in Section \ref{proofs} that the terms that depend
on the bandwidth $b$ are lower order terms in the expansion of $\hat
{\theta}_n-\theta.$ One can then minimise with respect to $b$ a
bound on these lower order terms. A minimiser, say $b^*,$ will
depend on the unknown true parameter~$\theta,$ also via
$x_{\theta}$ and~$x_{\theta}^{\prime},$ as well as on the error
variances $\sigma_1^2,\ldots,\sigma_d^2.$ However,
$\theta,x_{\theta},$ and~$x_{\theta}^{\prime}$ can be reestimated
via $\hat{\theta}_n,\hat{x},$ and $\hat{x}^{\prime}$ using a
different, pilot bandwidth $\tilde{b}.$ Of course, instead of
$\hat{x}$ and $\hat{x}^{\prime}$ the use of any other
nonparametric estimators of a regression function and its
derivative, for example, local polynomial estimators, see Section~1.6 of
\cite{tsybakov}, or the Gasser--M\"uller estimator, see
\cite{gasser}, is also a valid option. Error term variances can
be estimated via one of the methods described in \cite{hall} or
Section 5.6 of \cite{wasserman}. Once the pilot estimators of
$\theta,x_{\theta},$ and $x_{\theta}^{\prime}$ together with
estimators of $\sigma_1^2,\ldots,\sigma_d^2$ are available, these
can be plugged back into $b^{*}$ and in this way one obtains a
bandwidth $\hat{b}$ that estimates the optimal bandwidth $b^*.$
The final step would be computation of $\hat{\theta}_n$ with a~new
bandwidth $\hat{b}.$ Unfortunately, this method leads to
extremely cumbersome expressions and furthermore, since we are
minimising an upper bound on numerous remainder terms, it will
probably tend to oversmooth, that is, produce a bandwidth $b$ larger
than required. Moreover, the plug-in approach in general is
subject to some controversy having both supporters and critics,
see, for example, \cite{loader} and references therein. An
alternative to
the plug-in approach might be an approach based on one of the
resampling methods: cross-validation, jackknife, or bootstrap.
Computationally these resampling methods will be quite intensive.
Theoretical analysis of the properties of such bandwidth selectors
is a rather nontrivial task. Also a thorough simulation study is
needed before the practical value of different bandwidth selection
methods can be assessed. We do not address these issues here.

The next observation of this section concerns numerical
computation of our SME. The kernel-type
nonparametric regression estimates of $x_{\theta
j}, j=1,\ldots,d,$ can be quickly evaluated on any regular grid
of points $0\leq s_1\leq\cdots \leq s_m,$ for example, via techniques using
the Fast Fourier Transform (FFT) similar to those described in
Appendix D of \cite{wand2}. See also \cite{marron}. Furthermore, in
the match step of
the two-step estimation procedure the criterion function $M_{n,w}$
can be approximated by a finite sum by discretising the integral
in its definition. If $F$ is linear in $\theta_1,\ldots,\theta_p$
and is univariate,
then as already observed in~\cite{varah}, see pages 29 and 31,
cf. page 1262 in \cite{brunel} and page 1573 in \cite{liang},
this will lead to a~weighted linear least squares problem, which
can be solved in a routine fashion without using for example, random
search methods. This is a great simplification in comparison to
the ordinary least squares estimator, which moreover will still
tend to get trapped in local minima of the least squares criterion
function despite the fact that $F$ is linear in its parameters.

\begin{figure}[b]

\includegraphics{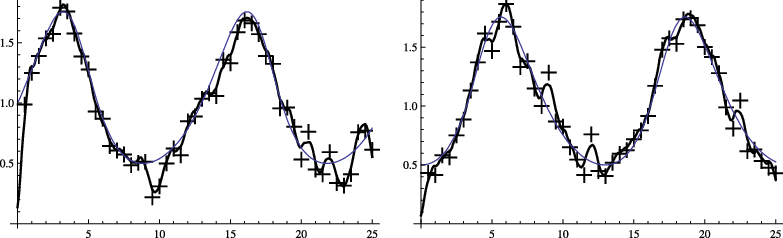}

\caption{Solution of the Lotka--Volterra system \protect\eqref{lotka} (thin
line) with parameter values $\theta_k=0.5,k=1,\ldots,4,$ and initial
condition $(x_1(0),x_2(0))=(1,0.5),$ observations $Y_{ij}$ given
by~\protect\eqref{lotkaobs} with $\epsilon_{ij}\sim N(0,0.01)$ (crosses) and the
estimates $\hat{x}_{j}$ computed with kernel \protect\eqref{kernelf}, weight
function \protect\eqref{weightf} and bandwidth $b=1.2$ (solid line). The left
panel corresponds to $x_{\theta1},$ the right to $x_{\theta2}.$}
\label{fig:sol}
\end{figure}

We conclude this section with two simple problems illustrating
parameter estimation for systems of ordinary differential equations via
the smooth and match method studied in the present paper. Our first
example deals with the Lotka--Volterra system that is a~basic model in
population dynamics. It describes evolution over time of the
populations of two species, predators and their preys. In mathematical
terms, the Lotka--Volterra model is described by a system consisting of
two ordinary differential equations and depending on the parameter
$\theta=(\theta_1,\theta_2,\theta_3,\theta_4)^{\mathrm{T}},$
%
\begin{equation}
\label{lotka}
\cases{
{ x}^{\prime}_1(t)=\theta_1 x_1(t)-\theta_2 x_1(t)x_2(t),\vspace*{2pt}\cr
{x}^{\prime}_2(t)=-\theta_3 x_2(t)+\theta_4 x_1(t)x_2(t).
}
\end{equation}
Here, $x_1$ represents the prey population and $x_2$ the predator
population. For additional information on the Lotka--Volterra system
see, for example, Section 6.2 in \cite{keshet}. We took $\theta
_k=0.5,k=1,\ldots,4$, and the initial condition
$(x_1(0),x_2(0))=(1,0.5).$ The solution to~\eqref{lotka} corresponding
to these parameter values is plotted in Figure \ref{fig:sol} with a
thin line. The left panel represents $x_{\theta1},$ the right panel
$x_{\theta2}.$
%
The solution components~$x_{\theta1}$ and~$x_{\theta2}$ are of
oscillatory nature and are out of phase of each other. Next, we
simulated a~small data set of size $n=50$ of observations on the
solution $x_{\theta}$ of \eqref{lotka} over the time interval
$[0,25]$ by taking an equidistant grid of time points $t_i=0.5 i$ for
$i=1,\ldots,50$ and setting
%
\begin{equation}
\label{lotkaobs}
Y_{ij}=x_{\theta j}(t_i)+\epsilon_{ij}, \quad i=1,\ldots,50,j=1,2,
\end{equation}
where the i.i.d. measurement errors $\epsilon_{ij}$ were generated
from the normal distribution $N(0,\sigma^2)$ with mean zero and
variance $\sigma^2=0.01.$ These observations $Y_{ij}$ are represented
by crosses in Figure \ref{fig:sol}.

\begin{figure}

\includegraphics{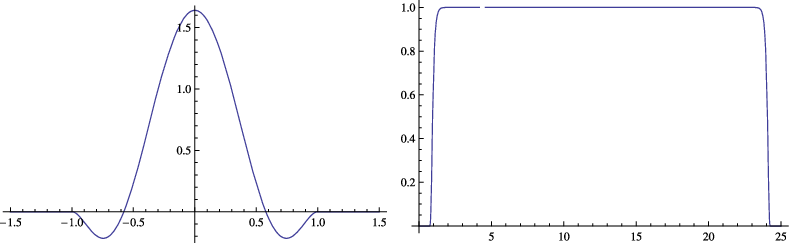}

\caption{Kernel $K$ from \protect\eqref{kernelf} (left panel) and weight
function $w$ from \protect\eqref{weightf} (right panel).}
\label{fig:kernel}
\vspace*{-3pt}
\end{figure}

The three required ingredients for the construction of an estimator
$\hat{\theta}_n$ are the kernel~$K,$ the weight function $w,$ and the
bandwidth $b.$ A general recipe for construction of kernels of an
arbitrary order $\alpha$ is given in Section 1.2.2 of \cite
{tsybakov} and is based on the use of polynomials that are orthonormal
in $L_2(-1,1)$ with weights. In particular, we used the ultraspherical
or Gegenbauer polynomials with weight function $v(t)=(1-t^2)^2
1_{[|t|\leq1]}$ and constructed the fourth order kernel with them.
Notice that our definition of the kernel of order $\alpha$ in
Condition \ref{cndkernel} is slightly different from the one in
Definition 1.3 of \cite{tsybakov}, cf. also the remark on page
6\vadjust{\goodbreak}
there. For ultraspherical polynomials, see Section 4.7 in \cite
{szego}. Our fourth order kernel took the form
%
\begin{equation}
\label{kernelf}
K(t)= \biggl( \frac{105}{64} - \frac{315}{64} t^2  \biggr) (1-t^2)^2
1_{[|t|\leq1]}.
\end{equation}
Notice that $K$ is a symmetric function. The kernel $K$ is plotted in
Figure \ref{fig:kernel} in the left panel.
%
An alternative here is to use the Gaussian-based kernels as in \cite
{schucany}, although they do not have a compact support. As far as the
weight function $w$ is concerned, any nonnegative function that is
equal to zero close to the end points of the interval $[0,25],$ is
equal to one on the greater part of the interval $[0,25]$ and is
smooth, could have been used. We opted to simply rescale and shift the function
\[
\lambda_{c,\beta}(t)=
\cases{
1, & \quad  if $|t|\leq c,$\cr
\exp[-\beta\exp[-\beta/(|t|-c)^2]/(|t|-1)^2], & \quad  if
$c<|t|<1,$\cr
0, & \quad  if $|t| \geq1,$
}\
\]
that arose in a different context in \cite{mcmurry}, see formula (3)
on page 552 there, so that it could have the required properties in
our context. We took the constants $c$ and $\beta$ to be equal to
$0.7$ and $0.5,$ respectively, and then set
%
\begin{equation}
\label{weightf}
w(t)=\lambda_{c,\beta} \biggl(1.05\frac{(t-12.5)}{12.5} \biggr).
\end{equation}
The function $w$ is plotted in the right panel of Figure \ref
{fig:kernel}. Finally, since in the present work construction of the
bandwidth selector is not our primary goal, we simply selected~$b$ by
hand and set it to 1.2.

\begin{figure}

\includegraphics{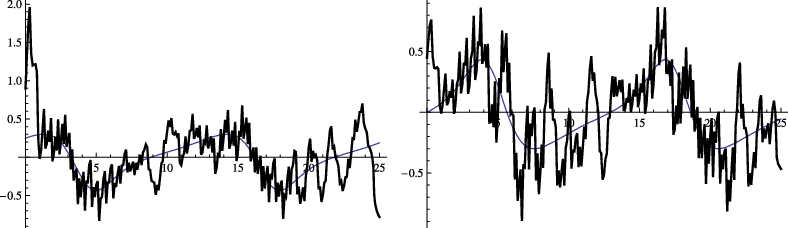}

\caption{Derivatives of the solution components $x_{\theta j}$ of the
Lotka--Volterra system \protect\eqref{lotka} (thin line) with parameter values
$\theta_k=0.5,k=1,\ldots,4,$ and initial condition
$(x_1(0),x_2(0))=(1,0.5),$ together with derivative estimates $\hat
{x}_{j}^{\prime}$ (solid line) computed with kernel \protect\eqref{kernelf},
weight function~\protect\eqref{weightf}, and bandwidth $b=1.2$ using
observations $Y_{ij}$ from \protect\eqref{lotkaobs}. The left panel
corresponds to~$\hat{x}_{1}^{\prime},$ the right panel to $\hat
{x}_{2}^{\prime}.$}
\label{fig:derivative}
\vspace*{6pt}
\end{figure}

The smooth and match estimation procedure was implemented in \emph
{Mathematica} 6.0, see \cite{mathematica}. We first evaluated the
kernel estimates of the regression functions $x_{\theta1}$ and
$x_{\theta2}$ at the equidistant grid of points $s_k=0.1 k$ with
$k=0,\ldots,249.$ With this number of grid points and the sample size
$n=50$ there was no need to use binning to compute the estimates and
moreover, binning would have probably resulted in a slower procedure,
cf. Figure 3(b) in \cite{marron}; so we did not employ it. However,
the fact that many of the kernel evaluations $K((s_k-t_i)/b)$ are
actually the same, cf. \cite{marron}, was taken into account and led
to savings in computation time above the naive implementation of the
Priestley--Chao estimator that would directly compute $K((t-t_i)/b).$
The estimates $\hat{x}_1$ and $\hat{x}_2$ are plotted in Figure \ref
{fig:sol} with a solid line, while the estimates $\hat{x}_1^{\prime}$
and $\hat{x}_2^{\prime}$ are plotted in Figure \ref{fig:derivative}.
%
Notice that the estimates $\hat{x}_1^{\prime}$ and $\hat
{x}_2^{\prime}$ are severely undersmoothed. We next approximated the
criterion function $M_{n,w}$ by a Riemann sum
\begin{eqnarray*}
&&\sum_{k=0}^{249} \bigl(\hat{x}_1^{\prime}(0.1 k) -\eta_1\hat{x}_1(0.1
k) +\eta_2 \hat{x}_1(0.1 k) \hat{x}_2(0.1 k)\bigr)^2 w(0.1k) 0.1 \\
&& \quad {}+
\sum_{k=0}^{249} \bigl( \hat{x}_2^{\prime}(0.1 k) +\eta_3 \hat{x}_2(0.1
k) - \eta_4 \hat{x}_1(0.1 k) \hat{x}_2(0.1 k) \bigr)^2 w(0.1 k)
0.1 .
\end{eqnarray*}
Note that when performing minimisation, the factor $0.1$ can be omitted
from both terms in the above display. The minimisation procedure
resulted in the estimate
\[
\hat{\theta}_n=(0.52,0.50,0.50,0.51)^{\mathrm{T}}.
\]
With our implementation, the total time needed for computation of the
estimate of $\theta$ (including time needed for kernel and weight
function evaluations, but excluding time needed for loading
observations) was about 0.5 seconds on a notebook with Intel(R)
Pentium(R) Dual CPU T3200 @ 2.00 GHz processor and 4.00 GB RAM. The
parameter estimates appear to be sufficiently accurate in this
particular case.\looseness=1

\begin{figure}

\includegraphics{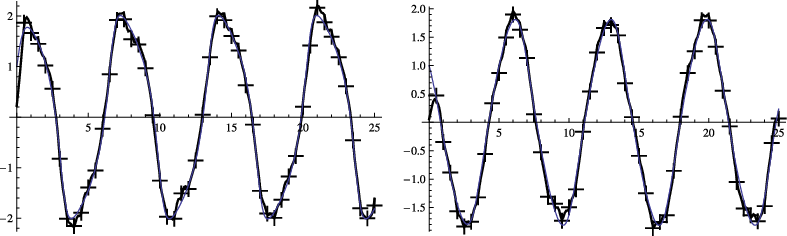}

\caption{Solution of the Van der Pol system \protect\eqref{pol} (thin line)
with parameter value $\theta=0.8$ and initial condition
$(x_1(0),x_2(0))=(1,1),$ observations $Y_{ij}$ given by \protect\eqref{polobs}
with $\epsilon_{ij}\sim N(0,0.01)$ (crosses) and the estimates $\hat
{x}_{j}$ computed with kernel \protect\eqref{kernelf}, weight function \protect\eqref
{weightf}, and bandwidth $b=1$ (solid line). The left panel corresponds
to $x_{\theta1}$ and the right to $x_{\theta2}.$}
\label{fig:pol}
\vspace*{-3pt}
\end{figure}

\begin{figure}

\includegraphics{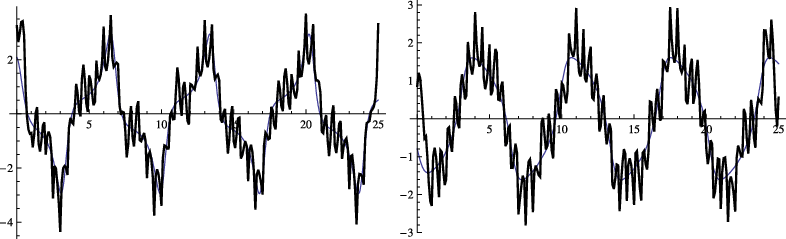}

\caption{Derivatives of the solution components $x_{\theta j}$ of the
Van der Pol system \protect\eqref{pol} (thin line) with parameter value
$\theta=0.8$ and initial condition $(x_1(0),x_2(0))=(1,1),$ together
with derivative estimates $\hat{x}_{j}^{\prime}$ (solid line)
computed with kernel \protect\eqref{kernelf}, weight function \protect\eqref
{weightf}, and bandwidth $b=1$ using observations $Y_{ij}$ from \protect\eqref
{polobs}. The left panel corresponds to $\hat{x}_{1}^{\prime}$ and
the right panel to $\hat{x}_{2}^{\prime}.$}
\label{fig:pol:der}
\vspace*{3pt}
\end{figure}

Our second example deals with the Van der Pol oscillator that describes
an electric circuit containing a nonlinear element, see page 333,
Problem 12 on page 365, and the references on page 373 in \cite
{keshet}. The corresponding system of ordinary differential equations
takes the form
%
\begin{equation}
\label{pol}
\cases{
{ x}^{\prime}_1(t)={\theta}^{-1} \bigl( x_1(t) - \frac{1}{3}
(x_1(t))^3 + x_2(t)  \bigr),\cr
{x}^{\prime}_2(t)=-\theta x_1(t).
}
\end{equation}
We took $\theta=0.8$ and the initial condition
$(x_1(0),x_2(0))=(1,1).$ The solution to~\eqref{pol} is of oscillatory
nature and the components $x_{\theta1}$ and $x_{\theta2}$ are out of
phase of each other. The solution is plotted in Figure \ref{fig:pol}
with a thin line.
%
We then simulated a~data set of size $n=50$ of observations on the
solution $x_{\theta}$ of \eqref{pol} over the time interval $[0,25]$
at an equidistant grid of time points $t_i=0.5 i,i=1,\ldots,50,$ by
setting
%
\begin{equation}
\label{polobs}
Y_{ij}=x_{\theta j}(t_i)+\epsilon_{ij}, \qquad i=1,\ldots,50,j=1,2,
\end{equation}
where the i.i.d. measurement errors $\epsilon_{ij}$ were generated
from the normal distribution $N(0,\sigma^2)$ with mean zero and
variance $\sigma^2=0.01.$ These observations $Y_{ij}$ are plotted with
crosses in Figure~\ref{fig:pol}. When computing the estimate $\hat
{\theta}_n,$ we used the same kernel and the same weight function as
in the previous example, while the bandwidth was set to $b=1.$ The
estimates of the solution components $x_{\theta1}$ and $x_{\theta2}$
are depicted by a solid line in Figure \ref{fig:pol}, while the
derivatives $x_{\theta1}^{\prime}$ and $x_{\theta2}^{\prime}$
together with their estimates are given in Figure \ref{fig:pol:der}.
%
The estimation procedure resulted in an estimate $\hat{\theta
}_n=0.83$ and the computation time was about $0.4$ seconds.

We intend to perform a more practically oriented study exploring
some of the ideas mentioned in this section in a separate
publication.

\section{Proofs}

When comparing two sequences $\alpha_n$ and $\beta_n$ of real
numbers, we will use the symbol $\lesssim,$ meaning $\alpha_n$ is
less or equal than $\beta_n$ up to a
universal multiplicative constant that is independent of index $n.$ The
symbol $\asymp$
will denote the fact that two sequences of real numbers are
asymptotically of the same order.\vadjust{\goodbreak}
\label{proofs}
\begin{pf*}{Proof of Proposition \ref{priestleyprop}}
We first prove \eqref{mup}. For any positive $\varepsilon$ by
Chebyshev's inequality we have
%
\begin{eqnarray}
\label{pr*}
P \Bigl(\sup_{t\in[\delta,1-\delta]}|\hat{\mu}_n(t)-\mu(t)| >
\varepsilon \Bigr)&\leq&\frac{2}{\varepsilon^2} \Bigl \{ \sup
_{t\in[\delta,1-\delta]}|\ex[\hat{\mu}_n(t)]-\mu(t)|^2\nonumber
\\
&&\hphantom{\frac{2}{\varepsilon^2} \Bigl \{}{} +\ex \Bigl[\sup_{t\in[\delta,1-\delta]}|\hat{\mu}_n(t)-\ex
[\hat{\mu}_n(t)]|^2 \Bigr]  \Bigr\}\\
&=&\frac{2}{\varepsilon^2}(T_1+T_2).
\nonumber
\end{eqnarray}
By formula \eqref{exapprox} from Appendix \ref{app1}, we can write
\[
\ex[\hat{\mu}_n(t)]-\mu(t)=\int_0^1 \mu(s) \frac1b
K \biggl(\frac{t-s}{b} \biggr)\,\mathrm{d}s-\mu(t)+\mathrm{O}\biggl(\frac{1}{nb^2} \biggr).
\]
For all $n$ large enough, we have $b\leq\delta,$ because
$b\rightarrow0.$ Then for all such $n,$ if
$t\in[\delta,1-\delta],$ a standard argument (cf. page 6 in
\cite{tsybakov}), namely Taylor's formula up to order $\alpha$
applied to $\mu$ and the moment conditions on the kernel $K$
formulated in Condition \ref{cndkernel}, yields
%
\begin{equation}
\label{7*}
\sup_{t\in[\delta,1-\delta]}|\ex[\hat{\mu}_n(t)]-\mu(t)| \leq
b^{\alpha} \frac{\| \mu^{(\alpha)} \|_{\infty}}{\alpha!}
\int_{-1}^{1} |u^{\alpha}K(u)|\,\mathrm{d}u + \mathrm{O}\biggl(\frac{1}{nb^2} \biggr).
\end{equation}
Next we turn to $T_2.$ With argumentation similar to that in the
proof of Theorem~1.8 of~\cite{tsybakov} and setting
\[
S_i(t)=\frac{t_i-t_{i-1}}{b} K \biggl(\frac{t-t_i}{b} \biggr),\qquad
N=n^2,  s_j=\frac{j}{N},
\]
for
$j=1,\ldots,N,$ we have
\begin{eqnarray*}
A&=&\sup_{t\in[\delta,1-\delta]}|\hat{\mu}_n(t)-\ex[\hat{\mu
}_n(t)]|\\
&=&\sup_{t\in[\delta,1-\delta]} \Biggl| \sum_{i=1}^n S_i(t) \epsilon
_i  \Biggr|\\
&\leq&\max_{1\leq j\leq N} \Biggl| \sum_{i=1}^n S_i(s_j) \epsilon_i
 \Biggr| +\sup_{t,t^{\prime}:|t-t^{\prime}|\leq N^{-1}} \Biggl|
\sum_{i=1}^n\bigl (S_i(t)-S_i(t^{\prime})\bigr) \epsilon_i  \Biggr|.
\end{eqnarray*}
By the mean value theorem and Condition \ref{cnd_times}, the inequality
\[
|S_i(t)-S_i(t^{\prime})|\lesssim\|K^{\prime}\|_{\infty}
\frac{1}{nb^2}|t-t^{\prime}|
\]
holds for any $t,t^{\prime}\in\mathbb{R},$ where
$\|K^{\prime}\|_{\infty}$ is finite. Hence, by the $c_2$-inequality
%
\begin{eqnarray}
\label{25half}
A^2&\leq& \Biggl( \max_{1\leq j\leq N}  \Biggl| \sum_{i=1}^n
\epsilon_i S_i(s_j)  \Biggr| +\sup_{t,t^{\prime}:|t-t^{\prime}|
\leq N^{-1}} \Biggl| \sum_{i=1}^n \bigl(S_i(t)-S_i(t^{\prime})\bigr) \epsilon_i
 \Biggr| \Biggr)^2\nonumber
 \\[-8pt]
 \\[-8pt]
&\lesssim&\max_{1\leq j\leq
N}|Z_j|^2+\frac{\|K^{\prime}\|_{\infty}^2}{n^2b^4N^2} \biggl(
\sum_{i=1}^n |\epsilon_i|  \biggr)^2,
\nonumber
\end{eqnarray}
where $Z_j=\sum_{i=1}^n S_i(s_j) \epsilon_i.$
Notice that
%
\begin{equation}
\label{25twth}
\frac{1}{n^2b^4N^2}\ex\Biggl [ \Biggl( \sum_{i=1}^n |\epsilon_i|
 \Biggr)^2 \Biggr]\leq
\frac{\ex[\epsilon_1^2]}{N^2b^4}=\frac{\sigma^2}{n^4b^4}=\mathrm{o}\Biggl(\frac{1}{nb} \Biggr).
\end{equation}
Moreover, we have
\begin{eqnarray*}
\ex[Z_j^2] & =&\sum_{i=1}^n \sigma^2(t_i-t_{i-1})^2 \biggl(\frac1b
K \biggl( \frac{t_i-s_j}{b}  \biggr) \biggr)^2\\
& \lesssim&\frac{\sigma^2\|K\|_{\infty}^2}{n^2b^2}\sum_{i=1}^n
1_{[|t_i-s_j|\leq b]}\\
& \leq&
\frac{1}{nb}c_1\sigma^2\|K\|_{\infty}^2\max\biggl (2,\frac
{1}{nb} \biggr),
\end{eqnarray*}
where the last inequality follows from Condition \ref{cnd_times}.
Since the $Z_j$'s, being a linear combination of independent Gaussian
random variables, are themselves Gaussian, Corollary~1.3 of~\cite{tsybakov}
and the fact that $N=n^2$ then entail
%
\begin{equation}
\label{25thf}
\ex \Bigl[ \max_{1\leq j\leq N}|Z_j|^2  \Bigr]=\mathrm{O}\biggl(\frac{\log
N}{nb} \biggr)=\mathrm{O}\biggl(\frac{\log n}{nb} \biggr).
\end{equation}
Combining \eqref{25half}, \eqref{25twth} and \eqref{25thf}, we obtain
%
\begin{equation}
\label{pr**}
\ex[A^2]=\mathrm{O}\biggl(\frac{\log n}{nb} \biggr).
\end{equation}
Taking
\[
\varepsilon= M  \Biggl( b^{\alpha} + \frac{1}{nb^2}
+\sqrt{\frac{\log n}{nb}} \Biggr)
\]
with an appropriate constant $M$ yields \eqref{mup} by
\eqref{pr*}, \eqref{7*}, and \eqref{pr**}.

As far as the proof of \eqref{muprimep} is concerned, it is very
much similar to the proof of \eqref{mup} and is therefore omitted.
This completes the proof of the proposition.
\end{pf*}

\begin{pf*}{Proof of Proposition \ref{consistency}}
From the definition of $M_{n,w}(\eta)$ and $M_w(\eta),$ the
elementary inequality
\[
\bigl|\|a_1\|^2 - \|a_2\|^2 \bigr|\leq\|a_1-a_2 \|(\|a_1\|+\|a_2\|)
\]
and the Cauchy--Schwarz inequality we have
%
\begin{eqnarray}
\label{brunel1}
&&\hspace*{-25pt}|M_{n,w}(\eta)-M_w(\eta)| \nonumber\\
&&\hspace*{-25pt} \quad \leq \biggl\{\int_{
0}^{1} \| \hat{x}^{\prime}(t) - F(x_{\theta}(t),\theta)
+ F({x}_{\theta}(t),\eta) - F(\hat{x}(t),\eta) \|^2w(t)\,\mathrm{d}t \biggr\}
^{1/2}\nonumber
\\[-12pt]
\\[-4pt]
&&\hspace*{-25pt} \qquad {}\times \Biggl\{ \sqrt{\int_{
0}^{1} \| \hat{x}^{\prime}(t) - F(\hat{x}(t),\eta)\|^2w(t)\,\mathrm{d}t}
+ \sqrt{ \int_0^1 \| F({x}_{\theta}(t),\theta) - F({x}_{\theta
}(t),\eta) \|^2w(t)\,\mathrm{d}t} \Biggr\}\nonumber\\
&&\hspace*{-25pt} \quad =\sqrt{T_1}\bigl(\sqrt{T_2}+\sqrt{T}_3\bigr).\nonumber
\end{eqnarray}
For $T_1$ we have that
%
\begin{eqnarray}
\label{star1}
&&T_1\leq2 \int_{
\delta}^{1-\delta} \| \hat{x}^{\prime}(t)   - F(x_{\theta
}(t),\theta)\|^2w(t)\,\mathrm{d}t\nonumber
\\[-8pt]
\\[-8pt]
&& \quad {}+2 \int_{ \delta}^{1-\delta} \| F({x}_{\theta}(t),\eta) -
F(\hat{x}(t),\eta) \|^2w(t)\,\mathrm{d}t.
\nonumber
\end{eqnarray}
By \eqref{xiprimep} it holds that
%
\begin{eqnarray}
\label{brunel1b}
&&\sup_{\eta\in\Theta}\int_{ \delta}^{1-\delta} \|
\hat{x}^{\prime}(t)   - F(x_{\theta}(t),\theta) \|^2
w(t)\,\mathrm{d}t \nonumber\\
&& \quad = \int_{ \delta}^{1-\delta} \| \hat{x}^{\prime}(t) -
x^{\prime}_{\theta}(t) \|^2
w(t)\,\mathrm{d}t\nonumber
\\[-8pt]
\\[-8pt]
&& \quad \leq\sum_{i=1}^d \sup_{t\in[\delta,1-\delta]}|\hat{x}^{\prime
}_i(t) -
x^{\prime}_{i,\theta}(t)|^2 \int_{ \delta}^{1-\delta} w(t)\,\mathrm{d}t\nonumber\\
&& \quad  \convp0.\nonumber
\end{eqnarray}
Moreover, by Lemma \ref{process} from Appendix \ref{app1} we obtain that
%
\begin{equation}
\label{4*} \sup_{\eta\in\Theta}\int_{\delta}^{1-\delta} \|
F(\hat{x}(t),\eta)-F({x}_{\theta}(t),\eta) \|^2w(t)\,\mathrm{d}t\convp0.
\end{equation}
Furthermore, $T_3=\mathrm{O}_P(1)$ as $n\rightarrow\infty,$ because
%
\begin{equation}
\label{spsi} \sup_{\eta\in\Theta}\int_{ \delta}^{1-\delta} \|
F({x}_{\theta}(t),\theta) - F({x}_{\theta}(t),\eta)
\|^2w(t)\,\mathrm{d}t<\infty
\end{equation}
by compactness of $\Theta$ and Condition \ref{cnd_solution}, and
$T_2=\mathrm{O}_P(1),$ because
%
\begin{equation}
\label{spsi2} \sup_{\eta\in\Theta}\int_{ \delta}^{1-\delta} \|
\hat{x}^{\prime}(t) - F(\hat{x}(t),\eta)\|^2w(t)\,\mathrm{d}t = \mathrm{O}_P(1)
\end{equation}
holds by the inequality
\begin{eqnarray*}
&&\int_{ \delta}^{1-\delta} \|
\hat{x}^{\prime}(t) - F(\hat{x}(t),\eta)\|^2w(t)\,\mathrm{d}t\\
&& \quad \lesssim
\int_{ \delta}^{1-\delta} \|
\hat{x}^{\prime}(t) - x_{\theta}^{\prime}(t)\|^2w(t)\,\mathrm{d}t
+
\int_{ \delta}^{1-\delta} \|
x_{\theta}^{\prime}(t)-F(x_{\theta}(t),\eta)\|^2w(t)\,\mathrm{d}t\\
&& \qquad {}+
\int_{ \delta}^{1-\delta} \|
F(x_{\theta}(t),\eta)-F(\hat{x}(t),\eta)\|^2w(t)\,\mathrm{d}t,
\end{eqnarray*}
Corollary \ref{corxi}, compactness of $\Theta,$ Condition \ref
{cnd_solution}, and Lemma \ref{process} from Appendix \ref{app1}.
Combination of \eqref{brunel1}--\eqref{spsi2} implies
that
\[
\sup_{\eta\in\Theta}|M_{n,w}(\eta)-M_w(\eta)|\convp0.
\]
The statement of the proposition then follows from this fact, the
identifiability condition~\eqref{identifiability}, and Theorem 5.7
of \cite{vdvaart} or more generally Corollary 3.2.3 in \cite{wellner}.
\end{pf*}

\begin{pf*}{Proof of Proposition \ref{functional}}
We interpret the derivative of a one-dimensional function of
$\theta$ as a row $p$-vector of partial derivatives and we denote
the $d \times p$-matrix of partial derivatives
${\partial}F_i(x,\theta)/{\partial\theta_j}, i=1,\ldots,d,
j=1,\ldots,p,$ by
$F^{\prime}_{\theta}(x,\theta).$

We have
\[
\frac{\mathrm{d}}{\mathrm{d}\theta} \| \hat{x}^{\prime}(t)-F(\hat{x}(t),\theta) \|^2
=-2\bigl(\hat{x}^{\prime}(t)-F(\hat{x}(t),\theta)\bigr)^{\mathrm{T}}F_{\theta}^{\prime
}(\hat{x}(t),\theta).
\]
With this in mind and interchanging the order of integration and
differentiation, we find that the derivative of $M_{n,w}$ from
\eqref{thetan} with respect to $\theta$ is given by
\[
-2\int_{\delta}^{1-\delta}\bigl(\hat{x}^{\prime}(t)-F(\hat
{x}(t),\theta)\bigr)^{\mathrm{T}}F_{\theta}^{\prime}(\hat{x}(t),\theta)w(t)\,\mathrm{d}t.
\]
Since $\theta$ is an interior point of $\Theta,$ there exists an
$\varepsilon>0,$ such that the open ball of radius~$\varepsilon$
around~$\theta$ is contained in $\Theta.$ Take
\[
G_n=\{ |\hat{\theta}_n-\theta|<\varepsilon/2 \}
\]
and notice that by consistency of $\hat{\theta}_n$ we have
$P(G_n)\rightarrow1$ as $n\rightarrow\infty.$ If $\hat{\theta}_n$
is a point of minimum of $M_{n,w},$ then necessarily
\[
1_{G_n}\int_{\delta}^{1-\delta}\bigl(\hat{x}^{\prime}(t)
-F(\hat{x}(t),\hat{\theta}_n)\bigr)^{\mathrm{T}}F_{\theta}^{\prime}(\hat
{x}(t),\hat{\theta}_n)w(t)\,\mathrm{d}t=0,
\]
where $0$ at the right-hand side denotes now a row $p$-vector with
all its entries equal to zero. The
latter display can be rearranged as
\begin{eqnarray*}
&&1_{G_n}\int_{\delta}^{1-\delta}
(F_{\theta}^{\prime}(\hat{x}(t),\hat{\theta}_n))^{\mathrm{T}} \times
\bigl\{\bigl(\hat{x}^{\prime}(t)-x^{\prime}_{\theta}(t)\bigr)+\bigl(F(x_{\theta
}(t),\theta)
- F(\hat{x}(t),\theta)\bigr)\\
&&\hphantom{1_{G_n}\int_{\delta}^{1-\delta}
(F_{\theta}^{\prime}(\hat{x}(t),\hat{\theta}_n))^{\mathrm{T}} \times
\bigl\{}{}+\bigl(F(\hat{x}(t),\theta)
-F(\hat{x}(t),\hat{\theta}_n)\bigr)\bigr\} w(t)\,\mathrm{d}t=0,
\end{eqnarray*}
where now $0$ on the right-hand side denotes a column $p$-vector
with its entries equal to zero. Note that we have
\[
F(\hat{x}(t),\theta) -F(\hat{x}(t),\hat{\theta}_n)=\int_0^1
F_{\theta}^{\prime}\bigl(\hat{x}(t),\hat{\theta}_n+\lambda
(\theta-\hat{\theta}_n)\bigr)\,\mathrm{d}\lambda (\theta-\hat{\theta}_n).
\]
Hence,
%
\begin{eqnarray}
\label{*}
&&1_{G_n}\int_{\delta}^{1-\delta}
(F_{\theta}^{\prime}(   \hat{x}(t),\hat{\theta}_n))^{\mathrm{T}} \int_0^1
F_{\theta}^{\prime}\bigl(\hat{x}(t),\hat{\theta}_n+\lambda
(\theta-\hat{\theta}_n)\bigr)\,\mathrm{d}\lambda  w(t)\,\mathrm{d}t  (\hat{\theta
}_n-\theta)\nonumber\\
&& \quad =1_{G_n}\int_{\delta}^{1-\delta}
(F_{\theta}^{\prime}(\hat{x}(t),\hat{\theta}_n))^{\mathrm{T}}
\bigl(\hat{x}^{\prime}(t)-x^{\prime}_{\theta}(t)\bigr)w(t)\,\mathrm{d}t\\
&& \qquad {}+1_{G_n}\int_{\delta}^{1-\delta}
(F_{\theta}^{\prime}(\hat{x}(t),\hat{\theta}_n))^{\mathrm{T}}
\bigl(F(x_{\theta}(t),\theta) - F(\hat{x}(t),\theta)\bigr) w(t)\,\mathrm{d}t\nonumber
\end{eqnarray}
holds. By the fact that $\hat{x}$ converges in probability as a
random element on $[\delta,1-\delta]$ to $x_{\theta},$ see
\eqref{xip}, consistency of $\hat{\theta}_n,$ continuity of
$F_{\theta}^{\prime},$ continuity of integration and the
continuous mapping theorem, see Theorem 18.11 in \cite{vdvaart},
we have
%
\begin{eqnarray}
\label{**}
&&\int_{\delta}^{1-\delta}
(F_{\theta}^{\prime}(\hat{x}(t),\hat{\theta}_n))^{\mathrm{T}} \int_0^1
 F_{\theta}^{\prime}\bigl(\hat{x}(t),\hat{\theta}_n+\lambda
(\theta-\hat{\theta}_n)\bigr)\,\mathrm{d}\lambda  w(t)\,\mathrm{d}t\nonumber
\\[-8pt]
\\[-8pt]
&& \quad \convp\int_{\delta}^{1-\delta}
(F_{\theta}^{\prime}(x_{\theta}(t),\theta))^{\mathrm{T}}
F_{\theta}^{\prime}(x_{\theta}(t),\theta) w(t)\,\mathrm{d}t = J_{\theta},
\nonumber
\end{eqnarray}
where $J_{\theta}$ is nonsingular by assumption
\eqref{brunel5}. Therefore, \eqref{*} shows that the asymptotic
behaviour of
$\hat{\theta}_n-\theta$ is given by
%
\begin{eqnarray}
\label{asymp}
&&J_{\theta}^{-1} \biggl(\int_{\delta}^{1-\delta}
(F_{\theta}^{\prime}(\hat{x}(t),\hat{\theta}_n))^{\mathrm{T}}
\bigl(\hat{x}^{\prime}(t)-x^{\prime}_{\theta}(t)\bigr)w(t)\,\mathrm{d}t\nonumber
\\[-8pt]
\\[-8pt]
&&\hphantom{J_{\theta}^{-1} \biggl(}{}+\int_{\delta}^{1-\delta}
(F_{\theta}^{\prime}(\hat{x}(t),\hat{\theta}_n))^{\mathrm{T}}
\bigl(F(x_{\theta}(t),\theta) - F(\hat{x}(t),\theta)\bigr) w(t)\,\mathrm{d}t \biggr).
\nonumber
\end{eqnarray}
It thus remains to be shown that this expression in fact reduces
to the right-hand side of \eqref{brunel3}. First of all, notice
that
%
\begin{eqnarray}
\label{6*}
&&\int_{\delta}^{1-\delta}
(F_{\theta}^{\prime} (\hat{x}(t),\hat{\theta}_n ))^{\mathrm{T}}
\bigl(\hat{x}^{\prime}(t)-x^{\prime}_{\theta}(t)\bigr)w(t)\,\mathrm{d}t\nonumber\\
&& \quad =\int_{\delta}^{1-\delta}
(F_{\theta}^{\prime}({x}_{\theta}(t),{\theta}))^{\mathrm{T}}
\bigl(\hat{x}^{\prime}(t)-x^{\prime}_{\theta}(t)\bigr)w(t)\,\mathrm{d}t\nonumber\\
&& \qquad {} +\int_{\delta}^{1-\delta}
\bigl(F_{\theta}^{\prime}(\hat{x}(t),\hat{\theta}_n)-F_{\theta
}^{\prime}({x}_{\theta}(t),{\theta})\bigr)^{\mathrm{T}}
\bigl(\hat{x}^{\prime}(t)-x^{\prime}_{\theta}(t)\bigr)w(t)\,\mathrm{d}t\\
&& \quad = - \int_{\delta}^{1-\delta}
 \biggl(\frac{\mathrm{d}}{\mathrm{d}t}[F_{\theta}^{\prime}({x}_{\theta}(t),{\theta
})w(t)] \biggr)^{\mathrm{T}}
\bigl(\hat{x}(t)-x_{\theta}(t)\bigr)\,\mathrm{d}t \nonumber\\
&& \qquad {} +\int_{\delta}^{1-\delta}
\bigl(F_{\theta}^{\prime}(\hat{x}(t),\hat{\theta}_n)-F_{\theta
}^{\prime}({x}_{\theta}(t),{\theta})\bigr)^{\mathrm{T}}
\bigl(\hat{x}^{\prime}(t)-x^{\prime}_{\theta}(t)\bigr)w(t)\,\mathrm{d}t,\nonumber
\end{eqnarray}
where the last equality follows by integration by parts and the
fact that $w(\delta)=w(1-\delta)=0.$ The first term at the right-hand
side of \eqref{6*}
appears also in the leading term
$\Gamma(\hat{x})-\Gamma(x_{\theta})$ of \eqref{brunel3}. We will
now show that the other term at the right-hand side of~\eqref{6*} is
negligible, that is,
\[
\int_{\delta}^{1-\delta}
\bigl(F_{\theta}^{\prime}(\hat{x}(t),\hat{\theta}_n)-F_{\theta
}^{\prime}({x}_{\theta}(t),{\theta})\bigr)^{\mathrm{T}}
\bigl(\hat{x}^{\prime}(t)-x^{\prime}_{\theta}(t)\bigr)w(t)\,\mathrm{d}t=\mathrm{o}_P(n^{-1/2}).
\]
By the Cauchy--Schwarz inequality,
\begin{eqnarray*}
 &&\biggl\|\int_{\delta}^{1-\delta}
\bigl(F_{\theta}^{\prime}(\hat{x}(t),\hat{\theta}_n)-F_{\theta
}^{\prime}({x}_{\theta}(t),{\theta})\bigr)^{\mathrm{T}}
\bigl(\hat{x}^{\prime}(t)-x^{\prime}_{\theta}(t)\bigr)w(t)\,\mathrm{d}t \biggr\|\\
&& \quad \leq \biggl\{ \int_{\delta}^{1-\delta} \|
F_{\theta}^{\prime}(\hat{x}(t),\hat{\theta}_n)-F_{\theta}^{\prime
}({x}_{\theta}(t),{\theta})\|^2
w(t)\,\mathrm{d}t  \biggr\}^{1/2}  \biggl\{ \int_{\delta}^{1-\delta} \|
\hat{x}^{\prime}(t)-x^{\prime}_{\theta}(t) \|^2 w(t)\,\mathrm{d}t
 \biggr\}^{1/2},
\end{eqnarray*}
where $\|\cdot\|$ denotes the Frobenius or the Hilbert--Schmidt
norm of a matrix (recall that it is submultiplicative). By
\eqref{xiprimep}, we have
\[
 \biggl\{ \int_{\delta}^{1-\delta} \|
\hat{x}^{\prime}(t)-x^{\prime}_{\theta}(t) \|^2
w(t)\,\mathrm{d}t \biggr\}^{1/2}=\mathrm{O}_P(1)\Biggl (
b^{\alpha-1}+\frac{1}{nb^3}+\sqrt{\frac{\log n}{nb^3}} \Biggr).
\]
Furthermore,
%
\begin{eqnarray}
\label{26}
&&\int_{\delta}^{1-\delta} \|
F_{\theta}^{\prime}(\hat{x}(t), \hat{\theta}_n)-F_{\theta
}^{\prime}({x}_{\theta}(t),{\theta})\|^2
w(t)\,\mathrm{d}t \nonumber\\
&& \quad \leq2\int_{\delta}^{1-\delta} \|
F_{\theta}^{\prime}(\hat{x}(t),\hat{\theta}_n)-F_{\theta}^{\prime
}({x}_{\theta}(t),\hat{\theta}_n)\|^2
w(t)\,\mathrm{d}t\nonumber
\\[-8pt]
\\[-8pt]
&& \qquad {}+2\int_{\delta}^{1-\delta} \|
F_{\theta}^{\prime}({x}_{\theta}(t),\hat{\theta}_n)-F_{\theta
}^{\prime}({x}_{\theta}(t),{\theta})\|^2
w(t)\,\mathrm{d}t\nonumber\\
&& \quad =2T_1+2T_2.
\nonumber
\end{eqnarray}
Denote $F_{\theta}^{\prime}(x,\theta)=A(x,\theta)=(a_{i,j}(x,\theta
))_{i,j}.$ For $T_1$, we have
\begin{eqnarray*}
T_1&=&\sum_{i,j}\int_{\delta}^{1-\delta} \bigl(a_{i,j}(\hat{x}(t),\hat
{\theta}_n)-a_{i,j}(x_{\theta}(t),\hat{\theta}_n)\bigr)^2w(t)\,\mathrm{d}t\\
&=&\sum_{i,j}\int_{\delta}^{1-\delta}  \biggl( \int_0^1 \frac
{\partial}{\partial x} a_{i,j}
\bigl(x_{\theta}(t)+\lambda\bigl(\hat{x}(t)-x_{\theta}(t)\bigr),\hat{\theta
}_n\bigr)\,\mathrm{d}\lambda  \bigl(\hat{x}(t)-x_{\theta}(t)\bigr)  \biggr)^2w(t)\,\mathrm{d}t\\
&\leq& \Bigl( \sup_{t\in[\delta,1-\delta]} \| \hat{x}(t)-x_{\theta
}(t) \|^2  \Bigr) \\
&&{}\times\sum_{i,j}\int_{\delta}^{1-\delta} \int_0^1  \biggl\|
\frac{\partial}{\partial x} a_{i,j}
\bigl(x_{\theta}(t)+\lambda\bigl(\hat{x}(t)-x_{\theta}(t)\bigr),\hat{\theta}_n\bigr)
 \biggr\|^2\,\mathrm{d}\lambda  w(t)\,\mathrm{d}t.
\end{eqnarray*}
By \eqref{xip}, as well as consistency of $\hat{\theta}_n,$
Condition \ref{cnd_solution} and the continuous mapping theorem,
the right-hand side in the last inequality is of order
\[
\mathrm{O}_P(1) \biggl\{ \biggl(b^{\alpha}+\frac{1}{nb^2} \biggr)^2+\frac{\log
n}{nb} \biggr\}.
\]
By a similar argument, the inequality
\begin{eqnarray*}
T_2&=&\int_{\delta}^{1-\delta} \|
F_{\theta}^{\prime}({x}_{\theta}(t),\hat{\theta}_n)
-F_{\theta}^{\prime}({x}_{\theta}(t),{\theta})\|^2
w(t)\,\mathrm{d}t\\
&\leq&\| \hat{\theta}_n-\theta\|^2
\sum_{i,j}\int_{\delta}^{1-\delta} \int_0^1  \biggl\|
\frac{\partial}{\partial\theta} a_{i,j}
\bigl(x_{\theta}(t),\theta+\lambda(\hat{\theta}_n-\theta)\bigr)  \biggr\|^2\,\mathrm{d}\lambda  w(t)\,\mathrm{d}t
\end{eqnarray*}
holds. Here with some natural abuse of notation we first
differentiate $a_{i,j}$ with respect to its second argument
$\theta$ and only afterwards evaluate the obtained derivative at
$x_{\theta}(t)$ and $\theta+\lambda(\hat{\theta}_n-\theta).$ Since
the integrals at the right-hand side of the above display are
bounded in probability, we then get
%
\begin{equation}
\label{5*}
 \biggl\{ \int_{\delta}^{1-\delta} \|
F_{\theta}^{\prime}({x}_{\theta}(t),\hat{\theta}_n)
-F_{\theta}^{\prime}({x}_{\theta}(t),{\theta})\|^2
w(t)\,\mathrm{d}t  \biggr\}^{1/2}=\mathrm{O}_P(\| \hat{\theta}_n-\theta\|).
\end{equation}
Now notice that \eqref{asymp} yields
\begin{eqnarray*}
&&\| \hat{\theta}_n-\theta\| \leq \mathrm{O}_P(1)
 \biggl( \biggl\|\int_{\delta}^{1-\delta}
(F_{\theta}^{\prime}(\hat{x}(t),\hat{\theta}_n))^{\mathrm{T}}
\bigl(\hat{x}^{\prime}(t)-x^{\prime}_{\theta}(t)\bigr)w(t)\,\mathrm{d}t \biggr\|\\
&& \hphantom{\| \hat{\theta}_n-\theta\| \leq \mathrm{O}_P(1)
 \biggl(}  {}+ \biggl\|\int_{\delta}^{1-\delta}
(F_{\theta}^{\prime}(\hat{x}(t),\hat{\theta}_n))^{\mathrm{T}}
\bigl(F({x}_{\theta}(t),\theta) -F(\hat{x}(t),{\theta})\bigr)
w(t)\,\mathrm{d}t \biggr\| \biggr).
\end{eqnarray*}
The Cauchy--Schwarz inequality then gives
\begin{eqnarray*}
\| \hat{\theta}_n-\theta\| &\leq& \mathrm{O}_P(1) \biggl \{
\int_{\delta}^{1-\delta} \|
F_{\theta}^{\prime}(\hat{x}(t),\hat{\theta}_n) \|^2w(t)\,\mathrm{d}t
 \biggr\}^{1/2}\\
&&{}\times \biggl\{ \int_{\delta}^{1-\delta}\|
\hat{x}^{\prime}(t)-x^{\prime}_{\theta}(t)
\|^2w(t)\,\mathrm{d}t \biggr\}^{1/2}\\
&&{}+ \mathrm{O}_P(1) \biggl \{ \int_{\delta}^{1-\delta} \|
F_{\theta}^{\prime}(\hat{x}(t),\hat{\theta}_n) \|^2w(t)\,\mathrm{d}t
 \biggr\}^{1/2}\\
&&\hphantom{{}+}{}\times\biggl \{ \int_{\delta}^{1-\delta} \| F({x}_{\theta
}(t),\theta)
-F(\hat{x}(t),{\theta}) \|^2w(t)\,\mathrm{d}t  \biggr\}^{1/2}.
\end{eqnarray*}
By a by now standard argument, that is,  \eqref{xip},
\eqref{xiprimep}, and the continuous mapping theorem, the right-hand
side can be further bounded to obtain
%
\begin{equation}
\label{star2} \| \hat{\theta}_n-\theta\| \leq \mathrm{O}_P(1)  \Biggl(
b^{\alpha-1}+\frac{1}{nb^3}+\sqrt{\frac{\log
n}{nb^3}} +
b^{\alpha}+\frac{1}{nb^2}+\sqrt{\frac{\log n }{nb}}
 \Biggr).
\end{equation}
Summarising the above results, we finally get that the second term
at the right-hand side of \eqref{6*} satisfies
\begin{eqnarray*}
&& \biggl\|\int_{\delta}^{1-\delta}
\bigl(F_{\theta}^{\prime}(\hat{x}(t),\hat{\theta}_n)-F_{\theta
}^{\prime}({x}_{\theta}(t),{\theta})\bigr)^{\mathrm{T}}
\bigl(\hat{x}^{\prime}(t)-x^{\prime}_{\theta}(t)\bigr)w(t)\,\mathrm{d}t \biggr\|\\
&& \quad \leq \mathrm{O}_P(1)
\Biggl (b^{\alpha-1}+\frac{1}{nb^3}+\sqrt{\frac{\log
n}{nb^3}} \Biggr)^2 =\mathrm{o}_P(n^{-1/2}),
\end{eqnarray*}
where the last equality follows from our
conditions on $b$. Here we also see that the condition $\alpha\geq
3$ is needed for the conclusion to hold.\vadjust{\goodbreak}

To conclude the proof, it
remains to consider the second term within brackets in \eqref{asymp}.
We have
%
\begin{eqnarray}
\label{asymp2}
&&\int_{\delta}^{1-\delta}
\bigl(F_{\theta}^{\prime}(\hat{x}(t),\hat{\theta}_n)\bigr)^{\mathrm{T}}
\bigl(F(x_{\theta}(t),\theta) - F(\hat{x}(t),\theta)\bigr) w(t)\,\mathrm{d}t\nonumber\\
&& \quad =\int_{\delta}^{1-\delta}
(F_{\theta}^{\prime}({x}_{\theta}(t),{\theta}))^{\mathrm{T}}
\bigl(F(x_{\theta}(t),\theta)
- F(\hat{x}(t),\theta)\bigr) w(t)\,\mathrm{d}t\\
&& \qquad {} +\int_{\delta}^{1-\delta}
\bigl(F_{\theta}^{\prime}(\hat{x}(t),\hat{\theta}_n)-F_{\theta
}^{\prime}({x}_{\theta}(t),{\theta})\bigr)^{\mathrm{T}}
\bigl(F(x_{\theta}(t),\theta) - F(\hat{x}(t),\theta)\bigr) w(t)\,\mathrm{d}t.\nonumber
\end{eqnarray}
This can be analysed in a by now routine fashion, but we provide
proofs. We first study the first term at the right-hand side. By a
standard argument, we have
\begin{eqnarray*}
&&\int_{\delta}^{1-\delta}
(F_{\theta}^{\prime}({x}_{\theta}(t),{\theta}))^{\mathrm{T}}
\bigl(F(x_{\theta}(t),\theta) - F(\hat{x}(t),\theta)\bigr) w(t)\,\mathrm{d}t\\
&& \quad =-\int_{\delta}^{1-\delta}
(F_{\theta}^{\prime}({x}_{\theta}(t),{\theta}))^{\mathrm{T}}\int_0^1
F_x^{\prime}\bigl(x_{\theta}(t)+\lambda
\bigl(\hat{x}(t)-x_{\theta}(t)\bigr),\theta\bigr)\,\mathrm{d}\lambda
  \bigl(\hat{x}(t)-x_{\theta}(t)\bigr)w(t)\,\mathrm{d}t\\
&& \quad =-\int_{\delta}^{1-\delta}
(F_{\theta}^{\prime}({x}_{\theta}(t),{\theta}))^{\mathrm{T}}
F_x^{\prime}(x_{\theta}(t),\theta)
\bigl(\hat{x}(t)-x_{\theta}(t)\bigr)w(t)\,\mathrm{d}t\\
&& \qquad {}-\int_{\delta}^{1-\delta}
(F_{\theta}^{\prime}({x}_{\theta}(t),{\theta}))^{\mathrm{T}}\int_0^1
\bigl[F_x^{\prime}\bigl(x_{\theta}(t)+\lambda\bigl(\hat{x}(t)-x_{\theta
}(t)\bigr),\theta\bigr)\\
&&\hphantom{={}-\int_{\delta}^{1-\delta}
(F_{\theta}^{\prime}({x}_{\theta}(t),{\theta}))^{\mathrm{T}}\int_0^1
\bigl[} \quad {}-F_x^{\prime}(x_{\theta}(t),\theta)\bigr]\,\mathrm{d}\lambda
  \bigl(\hat{x}(t)-x_{\theta}(t)\bigr)w(t)\,\mathrm{d}t\\
&& \quad =T_3+T_4.
\end{eqnarray*}
Recalling \eqref{brunel4}, we see that $T_3$ appears in the
leading term $\Gamma(\hat{x})-\Gamma(x_{\theta})$ in
\eqref{brunel3} and completes it together with the first term at
the right-hand side of \eqref{6*}. Next, we consider~$T_4.$
Introduce the notation
$F_x^{\prime}(x,\theta)=B(x,\theta)=(b_{i,j}(x,\theta))_{i,j}.$ We
have
\begin{eqnarray*}
 &&\biggl\| \int_0^1 \bigl[F_x^{\prime}\bigl(x_{\theta}(t)+\lambda(\hat{x}(t)
-x_{\theta}(t)),\theta\bigr)-F_x^{\prime}(x_{\theta}(t),\theta
)\bigr]\,\mathrm{d}\lambda
  \bigl(\hat{x}(t)-x_{\theta}(t)\bigr)  \biggr\|\\
&& \quad \leq \Bigl( \sup_{t\in[\delta,1-\delta]} \| \hat{x}(t)-x_{\theta
}(t) \|  \Bigr) \\
&& \qquad {}\times\int_0^1
\bigl\|F_x^{\prime}\bigl(x_{\theta}(t)+\lambda\bigl(\hat{x}(t)-x_{\theta
}(t)\bigr),\theta\bigr)
-F_x^{\prime}(x_{\theta}(t),\theta)\bigr\|\,\mathrm{d}\lambda\\
&& \quad \leq \Bigl( \sup_{t\in[\delta,1-\delta]} \| \hat{x}(t)-x_{\theta
}(t) \|  \Bigr)\\
&& \qquad {}\times\int_0^1 \sum_{i,j} \bigl|b_{i,j}\bigl(x_{\theta}(t)+\lambda\bigl(\hat{x}(t)
-x_{\theta}(t)\bigr),\theta\bigr) - b_{ij}(x_{\theta}(t),\theta)\bigr|\,\mathrm{d}\lambda\\
&& \quad \leq \Bigl( \sup_{t\in[\delta,1-\delta]} \| \hat{x}(t)-x_{\theta
}(t) \|  \Bigr)\\
&& \qquad {}\times\sum_{i,j} \int_0^1  \biggl\| \int_0^1 \frac{\partial
}{\partial x} b_{ij}
\bigl(x_{\theta}(t)+\kappa\lambda\bigl(\hat{x}(t)-x_{\theta}(t)\bigr),\theta
\bigr)\,\mathrm{d}\kappa\lambda\bigl(\hat{x}(t)-x_{\theta}(t)\bigr)  \biggr\|\,\mathrm{d}\lambda\\
&& \quad \leq \Bigl( \sup_{t\in[\delta,1-\delta]} \| \hat{x}(t)-x_{\theta
}(t) \|^2  \Bigr)\\
&& \qquad {}\times\sum_{i,j} \int_0^1 \int_0^1  \biggl\| \frac{\partial
}{\partial x}
b_{ij}\bigl(x_{\theta}(t)+\kappa\lambda\bigl(\hat{x}(t)-x_{\theta
}(t)\bigr),\theta\bigr) \biggr\|\,\mathrm{d}\kappa
\,\mathrm{d}\lambda,
\end{eqnarray*}
where in the last inequality we used the fact that
$0\leq\lambda\leq1.$ Since by convergence in probability of
$\hat{x}$ to $x_{\theta},$ Condition \ref{cnd_solution} and the
continuous mapping theorem the integrals on the right-hand side of
the above display are bounded in probability, it follows from~\eqref{xip} that $\|T_4\|$ is
\[
\mathrm{O}_P(1)\Biggl \{ \biggl (b^{\alpha}+\frac{1}{nb^2} \biggr)^2+\frac{\log
n}{nb} +  \biggl(b^{\alpha} + \frac{1}{nb^3}  \biggr) \sqrt{\frac
{\log n}{nb}}  \Biggr\}.
\]
This in turn is $\mathrm{o}_P(n^{-1/2})$ because
of the conditions on $b.$ Finally, we treat the second term at the
right-hand side of \eqref{asymp2}. By the Cauchy--Schwarz
inequality, its norm can be bounded by
\begin{eqnarray*}
&& \biggl\{\int_{\delta}^{1-\delta}
\|F_{\theta}^{\prime}(\hat{x}(t),\hat{\theta}_n)-F_{\theta
}^{\prime}({x}_{\theta}(t),{\theta})\|^2
w(t)\,\mathrm{d}t \biggr\}^{1/2}\\
&& \quad {}\times \biggl\{ \int_{\delta}^{1-\delta} \|F(x_{\theta}(t),\theta) -
F(\hat{x}(t),\theta)\|^2 w(t)\,\mathrm{d}t \biggr\}^{1/2}.
\end{eqnarray*}
Each of the terms at the right-hand side have already been treated
above, see \eqref{26} and~\eqref{star2}, and it follows that the
expression in the last display is $\mathrm{o}_P(n^{-1/2}).$ This
concludes the proof of Proposition \ref{functional}.
\end{pf*}

\begin{pf*}{Proof of Proposition \ref{linfunct}}
By a standard decomposition, we have
\begin{eqnarray*}
\ex\bigl[\bigl(\Delta(\hat{\mu}_n)-\Delta(\mu)\bigr)^2\bigr]&=&\bigl(\ex[\Delta(\hat{\mu
}_n)]-\Delta(\mu)\bigr)^2+\var[\Delta(\hat{\mu}_n)]\\
&=&T_1^2+T_2.
\end{eqnarray*}
The statement of the theorem will follow from Chebyshev's
inequality, provided we show that the right-hand side of the above
display is $\mathrm{O}( n^{-1}  ).$ For $T_1$, we have
\begin{eqnarray*}
|T_1|&=& \biggl| \int_{\mathbb{R}}
v(t)k(t)\bigl(\ex[\hat{\mu}_n(t)]-\mu(t)\bigr)\,\mathrm{d}t  \biggr|\\
&\leq&\sup_{t\in[\delta,1-\delta]}|\ex[\hat{\mu}_n(t)]-\mu
(t)|\int_{\mathbb{R}}
|v(t)k(t)|\,\mathrm{d}t\\
&=&\mathrm{O}\biggl(b^{\alpha}+\frac{1}{nb^2} \biggr),
\end{eqnarray*}
where the last equality follows from \eqref{7*}. Taking
$1/(2\alpha)\leq\gamma\leq1/4$ gives that $T_1$ is
$\mathrm{O}(n^{-1/2} ).$ We next consider $T_2.$ By independence
of the $\epsilon_i$'s, the fact that $\max_i
|t_i-t_{i-1}|\lesssim n^{-1},$ boundedness of $v$ and $k,$ and
integrability of $K,$ we have
\begin{eqnarray*}
T_2&=&\var \Biggl[ \sum_{i=1}^n (t_i-t_{i-1}) Y_i
\int_{\delta}^{1-\delta} v(t)k(t)\frac1b
K \biggl(\frac{t-t_i}{b} \biggr)\,\mathrm{d}t
 \Biggr]\\
&\lesssim&{\sigma^2} \sum_{i=1}^n
(t_i-t_{i-1})^2  \biggl(\int_{\delta}^{1-\delta}
v(t)k(t)\frac1b K\biggl (\frac{t-t_i}{b} \biggr)\,\mathrm{d}t \biggr)^2\\
&=&\mathrm{O}\biggl( \frac{1}{n}  \biggr).
\end{eqnarray*}
This
completes the proof of Proposition \ref{linfunct}.
\end{pf*}

\begin{pf*}{Proof of Theorem \ref{mainthm}}
The result is an easy consequence of Propositions \ref{functional}
and~\ref{linfunct}.
\end{pf*}

\begin{appendix}
\setcounter{equation}{0}
\section*{Appendix: Auxiliary results}
\label{apps}

\subsection{Technical lemmas}
\label{app1}

The proof of Proposition \ref{priestleyprop} is based on the
following two lemmas, which provide integral approximations to the
bias and variance of the estimator $\hat{\mu}_n$ and its
derivative $\hat{\mu}_n^{\prime}$ at a~point $t.$

\begin{lemma}
\label{lemmaintegral}
Let $\mu$ and $K$ be continuously differentiable and let $K$ be
supported on the interval $[-1,1].$ For any $t\in[0,1]$
%
\begin{equation}
\label{exapprox} \ex[\hat{\mu}_n(t)]=\int_0^1 \mu(s)\frac1b
K \biggl(\frac{t-s}{b} \biggr)\,\mathrm{d}s+\mathrm{O}\biggl(\frac{1}{nb^2} \biggr)
\end{equation}
holds in the regression model \eqref{regression}. The order bound on the
remainder term in \textup{\eqref{exapprox}} is uniform in $t\in[0,1].$
\end{lemma}
\begin{pf}
The proof is based on the Riemann sum approximation of the
integral. Since $\ex[\epsilon_i]=0,$ we have
\begin{eqnarray*}
\ex[\hat{\mu}_n(t)] &= & \int_0^1 \mu(s)\frac1b K \biggl(\frac
{t-s}{b} \biggr)\,\mathrm{d}s\\
&&{} - \int_0^1 \mu(s)\frac1b K \biggl(\frac{t-s}{b} \biggr)\,\mathrm{d}s +
\sum_{i=1}^n (t_i-t_{i-1}) \mu(t_i)\frac1b
K \biggl(\frac{t-t_i}{b} \biggr).
\end{eqnarray*}
The first term at the right-hand side of this expression is the
first term of \eqref{exapprox}. We will now establish an upper
bound on the difference of the other two terms. Using continuous
differentiability of $\mu$ and $K$ and the fact that
$\max_i|t_i-t_{i-1}|=\mathrm{O}(n^{-1}),$ we have
\begin{eqnarray*}
 &&\biggl| \int_0^1 \mu(s)\frac1b K \biggl(\frac{t-s}{b} \biggr)\,\mathrm{d}s
- \sum_{i=1}^n (t_i-t_{i-1}) \mu(t_i)\frac1b K \biggl(\frac
{t-t_i}{b} \biggr)  \biggr|\\[-2pt]
&& \quad =  \Biggl|\sum_{i=1}^n \int_{t_{i-1}}^{t_i}
 \biggl\{ \mu(s) \frac1b K \biggl(\frac{t-s}{b} \biggr) - \mu
(t_i)\frac1b K \biggl(\frac{t-t_i}{b} \biggr)  \biggr\}\,\mathrm{d}s  \Biggr|\\[-2pt]
&& \quad \leq\sum_{i=1}^n \int_{t_{i-1}}^{t_i}  \biggl| \mu(s)\frac1b K
 \biggl(\frac{t-s}{b} \biggr) - \mu(s)\frac1b K \biggl(\frac
{t-t_i}{b} \biggr)  \biggr|\,\mathrm{d}s\\[-2pt]
&& \qquad {}+ \sum_{i=1}^n \int_{t_{i-1}}^{t_i}  \biggl| \mu(s)
\frac1b K \biggl(\frac{t-t_i}{b} \biggr) - \mu(t_i)\frac1b K
\biggl(\frac{t-t_i}{b} \biggr)  \biggr|\,\mathrm{d}s\\[-2pt]
&& \quad \lesssim\frac{1}{nb^2}\| \mu\|_{\infty} \| K^{\prime} \|_{\infty} +
\frac{1}{nb} \| \mu^{\prime} \|_{\infty} \| K \|_{\infty},
\end{eqnarray*}
which is of order $n^{-1}b^{-2}.$ This establishes
\eqref{exapprox}.\vspace*{-2pt}
\end{pf}

The second lemma can be proved along the same lines as the
previous one and therefore we omit its proof. The existence of the
second derivative of $K$ is needed in the proof of this lemma.\vspace*{-2pt}
\begin{lemma}
\label{lemmaintegralder}
Let $\mu$ be continuously differentiable and let $K$ be twice
continuously differentiable and be supported on the interval
$[-1,1].$
For all $t\in[0,1]$
%
\begin{equation}
\label{exapproxder} \ex[\hat{\mu}_n^{\prime}(t)]=\int_0^1
\mu(s)\frac1{b^2}
K^{\prime} \biggl(\frac{t-s}{b} \biggr)\,\mathrm{d}s+\mathrm{O}\biggl(\frac
{1}{nb^3} \biggr)
\end{equation}
holds in the regression model \eqref{regression}. Furthermore, if
$b\leq\delta$ and
$t\in[\delta,1-\delta],$ then integration by parts yields
%
\begin{equation}
\label{exapproxder2} \ex[\hat{\mu}_n^{\prime}(t)]=\int_{-1}^1
\mu^{\prime}(t-bu)K(u)\,\mathrm{d}u+\mathrm{O}\biggl(\frac{1}{nb^3} \biggr).
\end{equation}
The order bounds on the remainder terms in \textup{\eqref{exapproxder}} and
\textup{\eqref{exapproxder2}} are uniform in $t$.\vspace*{-2pt}
\end{lemma}

The following lemma is used in the proof of Proposition
\ref{consistency}.\vspace*{-2pt}

\begin{lemma}
\label{process} Let the stochastic process
$X_n=(X_{n,\eta})_{\eta\in\Theta}$ be defined as
\[
X_n=(X_{n,\eta})_{\eta\in\Theta}= \biggl( \int_{\delta}^{1-\delta
} \|
F(\hat{x}(t),\eta)-F({x}_{\theta}(t),\eta) \|^2w(t)\,\mathrm{d}t
 \biggr)_{\eta\in\Theta}.
\]
Then under the conditions of Proposition \ref{consistency} we have
$X_n\convp0,$ where $0$ at the right-hand side denotes the zero
process on $\Theta$ and convergence\vadjust{\goodbreak} is understood as convergence for
random elements with values in the space $C(\Theta)$ of continuous
functions on $\Theta,$ which is equipped with the supremum norm.
\end{lemma}
\begin{pf}
To prove the lemma, we will verify the conditions of Theorem 18.14
of \cite{vdvaart}. By~\eqref{xip} and the continuous mapping
theorem, see Theorem 18.11 in \cite{vdvaart}, for every fixed~$\eta$ it holds that
%
\begin{equation}
\label{star3} \int_{\delta}^{1-\delta} \|
F(\hat{x}(t),\eta)-F({x}_{\theta}(t),\eta) \|^2w(t)\,\mathrm{d}t\convp0.
\end{equation}
Consequently, for any positive integer $k$ and any
$\eta_1,\ldots,\eta_k\in\Theta$ we have
\[
(X_{n,\eta_1},\ldots,X_{n,\eta_k})\rightsquigarrow
(\underbrace{0,\ldots,0}_k)
\]
and hence condition (i) of Theorem 18.14 in \cite{vdvaart} is
satisfied. Introduce
\[
G=\bigcap_{j=1}^d  \Bigl\{ \sup_{t\in[\delta,1-\delta]}|
\hat{x}_j(t)-x_{\theta j}(t) | \leq\beta \Bigr\}
\]
and notice
\[
G^c=\bigcup_{j=1}^d  \Bigl\{ \sup_{t\in[\delta,1-\delta]}|
\hat{x}_j(t)-x_{\theta j}(t) | > \beta \Bigr\}.
\]
For any positive $\varepsilon$ and $\beta$ and any partition
$\Theta_1,\ldots,\Theta_{m}$ of $\Theta$, we have
%
\begin{eqnarray}
\label{number1}
&&P\Bigl ( \sup_{\ell} \sup_{\eta,\zeta\in\Theta_{\ell}}
|X_{n,\eta}-X_{n,\zeta}| \geq\varepsilon \Bigr)\nonumber
\\[-8pt]
\\[-8pt]
&& \quad \leq P \Bigl( \sup_{\ell}\sup_{\eta,\zeta\in\Theta_{\ell}}
|X_{n,\eta}-X_{n,\zeta}| \geq\varepsilon; G  \Bigr)
+P (G^c ).
\nonumber
\end{eqnarray}
By \eqref{xip}, we know that
%
\begin{equation}
\label{number2} \lim_{n\rightarrow\infty} P (G^c ) \leq
\lim_{n\rightarrow\infty} \sum_{j=1}^d P \Bigl(
\sup_{t\in[\delta,1-\delta]}| \hat{x}_j(t)-x_{\theta j}(t) | >
\beta \Bigr) = 0.
\end{equation}
We will now show that for arbitrarily small positive $\rho$ and
$\varepsilon$ there exists a partition
$\Theta_1,\ldots,\Theta_{m}$ of $\Theta,$ such that
\[
\limsup_{n\rightarrow\infty}P \Bigl(
\sup_{\ell}\sup_{\eta,\zeta\in\Theta_{\ell}}
|X_{n,\eta}-X_{n,\zeta}| \geq\varepsilon; G  \Bigr)\leq\rho.
\]
Together with \eqref{number1} and \eqref{number2} this will imply
condition (ii) of Theorem 18.14 in \cite{vdvaart} and hence also
the fact that $X_n$ converges weakly to zero. The statement of the lemma will then be a
simple consequence of the fact that convergence to a constant in
distribution and in probability are equivalent, see Theorem 18.10 of \cite{vdvaart}.\vadjust{\goodbreak}

Notice that
\begin{eqnarray*}
&&|X_{n,\eta}-X_{n,\zeta}|\\[-2pt]
&& \quad \leq\int_{\delta}^{1-\delta} \|
F(\hat{x}(t),\eta)-F({x}_{\theta}(t),\eta) - F(\hat{x}(t),\zeta
)+F({x}_{\theta}(t),\zeta) \|\\[-2pt]
&&\hphantom{\int_{\delta}^{1-\delta}} \qquad {}\times\bigl(\| F(\hat{x}(t),\eta)-F({x}_{\theta}(t),\eta) \| + \|
F(\hat{x}(t),\zeta)-F({x}_{\theta}(t),\zeta) \|\bigr)w(t)\,\mathrm{d}t\\[-2pt]
&& \quad \leq \biggl\{ \int_{\delta}^{1-\delta} \|
F(\hat{x}(t),\eta)-F({x}_{\theta}(t),\eta) - F(\hat{x}(t),\zeta
)+F({x}_{\theta}(t),\zeta) \|^2w(t)\,\mathrm{d}t  \biggr\}^{1/2}\\[-2pt]
&& \qquad {}\times \biggl\{ \int_{\delta}^{1-\delta}\bigl (\| F(\hat{x}(t),\eta)
-F({x}_{\theta}(t),\eta) \| + \| F(\hat{x}(t),\zeta)-F({x}_{\theta
}(t),\zeta) \|\bigr)^2w(t)\,\mathrm{d}t  \biggr\}^{1/2}\\[-2pt]
&& \quad = \sqrt{T_3}\sqrt{T_4}.
\end{eqnarray*}
For $T_3$, we have
\begin{eqnarray*}
T_3&\leq&2 \int_{\delta}^{1-\delta} \|
F(\hat{x}(t),\eta) - F(\hat{x}(t),\zeta)\|^2w(t)\,\mathrm{d}t\\[-2pt]
&&{}+2 \int_{\delta}^{1-\delta} \| F({x}_{\theta}(t),\eta) -
F({x}_{\theta}(t),\zeta) \|^2w(t)\,\mathrm{d}t.
\end{eqnarray*}
Restricting $\omega$'s
from the sample space $\Omega$ to the set $G,$ we get
\begin{eqnarray*}
T_3 &\leq&2 \int_{\delta}^{1-\delta} \int_0^1 \bigl\| F_{\theta}^{\prime}
\bigl(\hat{x}(t),\zeta+\lambda(\eta-\zeta)\bigr) \bigr\|^2\,\mathrm{d}\lambda  \|\eta
-\zeta\|^2 w(t)\,\mathrm{d}t\\[-2pt]
&&{}+2 \int_{\delta}^{1-\delta} \int_0^1 \bigl\| F_{\theta}^{\prime
}\bigl({x}_{\theta}(t),\zeta+\lambda(\eta-\zeta)\bigr)
\bigr\|^2\,\mathrm{d}\lambda  \|\eta-\zeta\|^2 w(t)\,\mathrm{d}t\\[-2pt]
&\leq&4 \|\eta-\zeta\|^2 \int_{\delta}^{1-\delta} w(t)\,\mathrm{d}t
\mathop{\mathop{\sup}_{\|x_j\| \leq\| x_{\theta j} \|_{\infty} +
\beta,j=1,\ldots,d
}}_{\nu\in\Theta} \| F_{\theta}^{\prime}(x,\nu) \| =
C(\beta,w,\theta,\Theta) \|\eta-\zeta\|^2
\end{eqnarray*}
on the set $G.$ Notice that $C(\beta,w,\theta,\Theta)$ is a finite
constant, because $\|F_{\theta}^{\prime}(x,\nu)\|$ is continuous
and its supremum is taken over a compact set. By similar
techniques, one can show that on the set $G$ one has $T_4\leq
C^{\prime}(\beta,w,\theta,\Theta)$ for some constant
$C^{\prime}(\beta,w,\theta,\Theta),$ which depends only on
$\beta,w,\theta,$ and $\Theta.$ Consequently,
%
\begin{eqnarray}
\label{prob}
&&P \Bigl( \sup_{\ell} \sup_{\eta,\zeta\in\Theta_{\ell}}
|X_{n,\eta}-X_{n,\zeta}| \geq\varepsilon; G  \Bigr)\nonumber
\\[-9pt]
\\[-9pt]
&& \quad \leq P  \Bigl( \sup_{\ell}\sup_{\eta,\zeta\in\Theta_{\ell}}
{\sqrt{C(\beta,w,\theta,\Theta)C^{\prime}(\beta,w,\theta,\Theta
)}}
\| \eta- \zeta\| \geq\varepsilon \Bigr).
\nonumber
\end{eqnarray}
Now take a partition $\Theta_1,\ldots,\Theta_m$ of $\Theta,$ such
that for all $\ell=1,\ldots,m$
\[
0<\operatorname{diam}\Theta_{\ell}<\frac{\varepsilon}
{\sqrt{C(\beta,w,\theta,\Theta)C^{\prime}(\beta,w,\theta,\Theta)}}\vadjust{\goodbreak}
\]
holds, where $\operatorname{diam}\Theta_{\ell}$ denotes the
diameter of the set $\Theta_{\ell}.$ Observe that since
$\Theta\subset{\mathbb{R}}^p$ is compact, there indeed exists a
finite $m$ for which this is satisfied. The right-hand side of~\eqref{prob}
for such a partition is zero and consequently the
conditions (i) and (ii) of Theorem~18.14 of \cite{vdvaart} hold.
This completes the proof of the lemma.
\end{pf}

\subsection{Bounded measurement errors}
\label{app2}

Here we state and prove a modification of Proposition
\ref{priestleyprop} for the case when the $\epsilon_i$'s are
bounded.

\begin{prop}
\label{priestleypropbounded} In the regression model
\eqref{regression}, replace the assumption of Gaussianity of the
$\epsilon_i$'s by $|\epsilon_i|\leq C$ for some constant $C>0$ and
suppose Condition \ref{cndkernel} holds.
\begin{longlist}[(ii)]
\item[(i)] If $\mu$ is $\alpha\geq1$ times continuously differentiable
and $b\rightarrow0$ as $n\rightarrow\infty,$ then
%
\begin{equation}
\label{mupbounded}
\sup_{t\in[\delta,1-\delta]}|\hat{\mu}_n(t)-\mu(t)| = \mathrm{O}_P\Biggl(
b^{\alpha}+\frac{1}{nb^2}+\sqrt{\frac{\log n}{nb}}
 \Biggr).
\end{equation}

\item[(ii)] If $\mu$ is $\alpha\geq2$ times continuously differentiable
and $b\rightarrow0$ as $n\rightarrow\infty,$ then
%
\begin{equation}
\label{muprimepbounded}
\sup_{t\in[\delta,1-\delta]}|\hat{\mu}_n^{\prime}(t)-\mu
^{\prime}(t)|
= \mathrm{O}_P\Biggl(
b^{\alpha-1}+\frac{1}{nb^3}+\sqrt{\frac{\log
n}{nb^3}}  \Biggr)
\end{equation}
is valid. Moreover, $\hat{\mu}_n$ and $\hat{\mu}_n^{\prime}$ are
consistent on $[\delta,1-\delta],$ if $nb^3/\log
n\rightarrow\infty$ holds additionally.
\end{longlist}
\end{prop}

\begin{pf}
The proof of \eqref{mupbounded} follows the same steps as the
proof of \eqref{mup}. The only difference is that we need to show
that
%
\begin{equation}
\label{logbound} \ex \Bigl[\max_{1\leq j \leq N} |Z_j|^2 \Bigr]=
\mathrm{O}\biggl(\frac{\log n}{nb} \biggr)
\end{equation}
holds also for bounded $\epsilon_i$'s and not only for the
Gaussian $\epsilon_i$'s. To this end, we will use some results from
Chapter 2.2 of \cite{wellner}. Let $\eta$ be a nondecreasing and
convex function on $[0,\infty),$ such that $\eta(0)=0.$ The Orlicz
norm $\|X\|_{\eta}$ of a random variable $X$ is defined as
\[
\|X\|_{\eta}=\inf \biggl\{ C>0: \ex \biggl[ \eta \biggl( \frac{|X|}{C}
 \biggr)  \biggr]\leq1  \biggr\}.
\]
A particular $\eta$ that we will use is $\eta(x)=\exp(x^2)-1.$ Since
the $\epsilon_i$'s have mean zero and are
bounded, for any $x>0$ Hoeffding's inequality, see Theorem 2 in
\cite{hoeffding}, implies
\[
P(|Z_j|>x)\leq2 \exp \Biggl(-2 x^2\Big/ \Biggl(\sum_{i=1}^n C^2
(S_i(s_j))^2  \Biggr) \Biggr).\vadjust{\goodbreak}
\]
By Condition \ref{cnd_times}
\begin{eqnarray*}
C^2 \sum_{i=1}^n (S_i(s_j))^2 & \lesssim & C^2 \|K\|_{\infty}^2
\frac{1}{n^2 b^2}\sum_{i=1}^n
1_{[|s_j-t_i|\leq b]}\\[-2pt]
& \leq &\frac{1}{nb} C^2 \|K\|_{\infty}^2 c_1
\max \biggl(2,\max_n\frac{1}{nb} \biggr)= \frac{1}{C_0 nb}
\end{eqnarray*}
holds. Thus, the inequality
\[
P(|Z_j|>x)\leq2 \exp(-2C_0 nb x^2)
\]
is valid. By Lemma 2.2.1 of \cite{wellner}, it then follows that
%
\begin{equation}
\label{display*} \max_j\| Z_j \|_{\eta}\leq
\frac{C_1}{\sqrt{nb}},
\end{equation}
where $C_1$ depends on $C_0$ only. Let $\|X\|_2$ denote the $L_2$
norm of a random variable $X,$ that is, $\|X\|_2=\sqrt{\ex[X^2]}.$
Notice that the inequality
%
\begin{equation}
\label{orlicz} \|X\|_2\leq \|X\|_{\eta},
\end{equation}
holds, because of $\eta(x)\geq x^2.$ The inequalities
\eqref{display*} and \eqref{orlicz} combined with Lemma~2.2.2 of
\cite{wellner} yield that
\[
\sqrt{\ex \Bigl[\max_{1\leq j \leq N} |Z_j|^2 \Bigr]} \leq
\frac{C_3}{\sqrt{nb}}  \eta^{-1}(N),
\]
where the constant $C_3$ is independent of $N.$ Now notice that
for $N\geq4$
\[
\eta^{-1}(N) = \sqrt{\log(N+1)} \leq\sqrt{\log(N^2)} = 2
\sqrt{\log n}.
\]
Hence, \eqref{logbound} holds and this completes the proof of
\eqref{mupbounded}. Formula \eqref{muprimepbounded} can be proved
in a similar fashion.\vspace*{-2pt}
\end{pf}
\end{appendix}

\section*{Acknowledgements}\vspace*{-2pt}
The research reported here was started when the first author was a
Postdoc at EURANDOM, Eindhoven,
and the second one was a Senior Fellow there. The authors would like to
thank the referees and the Associate Editor for their suggestions and
comments on an earlier draft of the paper.\vspace*{-2pt}

%

\printhistory


\begin{thebibliography}{59}

\bibitem{arnold}
%
\begin{bbook}[mr]
\bauthor{\bsnm{Arnol'd},~\bfnm{V.~I.}\binits{V.I.}}
(\byear{1973}).
\btitle{Ordinary Differential Equations}.
\baddress{Cambridge, MA}: \bpublisher{MIT Press}.
\bid{mr={0361233}}
\end{bbook}
%
\endbibitem

\bibitem{bellman}
%
\begin{barticle}[mr]
\bauthor{\bsnm{Bellman},~\bfnm{Richard}\binits{R.}} \AND
\bauthor{\bsnm{Roth},~\bfnm{Robert~S.}\binits{R.S.}}
(\byear{1971}).
\btitle{The use of splines with unknown end points in the
identification of
systems}.
\bjournal{J. Math. Anal. Appl.}
\bvolume{34}
\bpages{26--33}.
\bid{issn={0022-247X}, mr={0277269}}
\end{barticle}\vadjust{\goodbreak}
%
\endbibitem

\bibitem{benedetti}
%
\begin{barticle}[mr]
\bauthor{\bsnm{Benedetti},~\bfnm{Jacqueline~K.}\binits{J.K.}}
(\byear{1977}).
\btitle{On the nonparametric estimation of regression functions}.
\bjournal{J. Roy. Statist. Soc. Ser. B}
\bvolume{39}
\bpages{248--253}.
\bid{issn={0035-9246}, mr={0494656}}
\end{barticle}
%
\endbibitem

\bibitem{bkrw}
%
\begin{bbook}[mr]
\bauthor{\bsnm{Bickel},~\bfnm{Peter~J.}\binits{P.J.}},
\bauthor{\bsnm{Klaassen},~\bfnm{Chris A.~J.}\binits{C.A.J.}},
\bauthor{\bsnm{Ritov},~\bfnm{Ya'acov}\binits{Y.}} \AND
\bauthor{\bsnm{Wellner},~\bfnm{John~A.}\binits{J.A.}}
(\byear{1998}).
\btitle{Efficient and Adaptive Estimation for Semiparametric Models}.
\baddress{New York}: \bpublisher{Springer}.
\bid{mr={1623559}}
\end{bbook}
%
\endbibitem

\bibitem{bickel}
%
\begin{barticle}[mr]
\bauthor{\bsnm{Bickel},~\bfnm{Peter~J.}\binits{P.J.}} \AND
\bauthor{\bsnm{Ritov},~\bfnm{Ya'acov}\binits{Y.}}
(\byear{2003}).
\btitle{Nonparametric estimators which can be ``plugged-in''}.
\bjournal{Ann. Statist.}
\bvolume{31}
\bpages{1033--1053}.
\bid{doi={10.1214/aos/1059655904}, issn={0090-5364}, mr={2001641}}
\end{barticle}
%
\endbibitem

\bibitem{bock}
%
\begin{bincollection}[mr]
\bauthor{\bsnm{Bock},~\bfnm{H.~G.}\binits{H.G.}}
(\byear{1983}).
\btitle{Recent advances in parameter identification techniques for {ODE}}.
In \bbooktitle{Numerical Treatment of Inverse Problems in Differential and
Integral Equations ({H}eidelberg, 1982)}.
\bseries{Progr. Sci. Comput.}
\bvolume{2}
\bpages{95--121}.
\baddress{Boston, MA}: \bpublisher{Birkh\"auser}.
\bid{mr={0714563}}
\end{bincollection}
%
\endbibitem

\bibitem{brunel}
%
\begin{barticle}[mr]
\bauthor{\bsnm{Brunel},~\bfnm{Nicolas J-B.}\binits{N.J.B.}}
(\byear{2008}).
\btitle{Parameter estimation of {ODE}'s via nonparametric estimators}.
\bjournal{Electron. J. Stat.}
\bvolume{2}
\bpages{1242--1267}.
\bid{doi={10.1214/07-EJS132}, issn={1935-7524}, mr={2471285}}
\end{barticle}
%
\endbibitem

\bibitem{chou}
%
\begin{barticle}[mr]
\bauthor{\bsnm{Chou},~\bfnm{I-Chun}\binits{I.C.}} \AND
\bauthor{\bsnm{Voit},~\bfnm{Eberhard~O.}\binits{E.O.}}
(\byear{2009}).
\btitle{Recent developments in parameter estimation and structure
identification of biochemical and genomic systems}.
\bjournal{Math. Biosci.}
\bvolume{219}
\bpages{57--83}.
\bid{doi={10.1016/j.mbs.2009.03.002}, issn={0025-5564}, mr={2537454}}
\end{barticle}
%
\endbibitem

\bibitem{keshet}
%
\begin{bbook}[mr]
\bauthor{\bsnm{Edelstein-Keshet},~\bfnm{Leah}\binits{L.}}
(\byear{2005}).
\btitle{Mathematical Models in Biology}.
\bseries{Classics in Applied Mathematics}
\bvolume{46}.
\baddress{Philadelphia, PA}: \bpublisher{SIAM}.
\bid{mr={2131632}}
\end{bbook}
%
\endbibitem

\bibitem{ellner}
%
\begin{barticle}[auto:STB|2011/09/12|07:03:23]
\bauthor{\bsnm{Ellner},~\bfnm{S.~P.}\binits{S.P.}},
\bauthor{\bsnm{Seifu},~\bfnm{Y.}\binits{Y.}} \AND
\bauthor{\bsnm{Smith},~\bfnm{R.~H.}\binits{R.H.}}
(\byear{2002}).
\btitle{Fitting population dynamic models to time-series data by gradient
matching}.
\bjournal{Ecology}
\bvolume{83}
\bpages{2256--2270}.
\end{barticle}
%
\endbibitem

\bibitem{marron}
%
\begin{barticle}[auto:STB|2011/09/12|07:03:23]
\bauthor{\bsnm{Fan},~\bfnm{J.}\binits{J.}} \AND
\bauthor{\bsnm{Marron},~\bfnm{J.~S.}\binits{J.S.}}
(\byear{1994}).
\btitle{Fast implementations of nonparametric curve estimators}.
\bjournal{J. Comput. Graph. Stat.}
\bvolume{3}
\bpages{35--56}.
\end{barticle}
%
\endbibitem

\bibitem{feinberg}
%
\begin{bmisc}[auto:STB|2011/09/12|07:03:23]
\bauthor{\bsnm{Feinberg},~\bfnm{M.}\binits{M.}}
(\byear{1979}).
\bhowpublished{Lectures on chemical reaction networks. Lectures delivered at the
Mathematics Research Center.
 Univ. Wisconsin-Madison. Available at
 \texttt{\href{http://www.che.eng.ohio-state.edu/\textasciitilde feinberg/LecturesOnReactionNetworks}
 {http://}
 \href{http://www.che.eng.ohio-state.edu/\textasciitilde feinberg/LecturesOnReactionNetworks}
 {www.che.eng.ohio-state.edu/\textasciitilde feinberg/LecturesOnReactionNetworks}}}.
\end{bmisc}
%
\endbibitem
\mbox{}

\bibitem{gasser}
%
\begin{barticle}[mr]
\bauthor{\bsnm{Gasser},~\bfnm{Theo}\binits{T.}} \AND
\bauthor{\bsnm{M{\"u}ller},~\bfnm{Hans-Georg}\binits{H.G.}}
(\byear{1984}).
\btitle{Estimating regression functions and their derivatives by the kernel
method}.
\bjournal{Scand. J. Statist.}
\bvolume{11}
\bpages{171--185}.
\bid{issn={0303-6898}, mr={0767241}}
\end{barticle}
%
\endbibitem

\bibitem{gasser2}
%
\begin{barticle}[mr]
\bauthor{\bsnm{Gasser},~\bfnm{T.}\binits{T.}},
\bauthor{\bsnm{M{\"u}ller},~\bfnm{H.~G.}\binits{H.G.}} \AND
\bauthor{\bsnm{Mammitzsch},~\bfnm{V.}\binits{V.}}
(\byear{1985}).
\btitle{Kernels for nonparametric curve estimation}.
\bjournal{J.~Roy. Statist. Soc. Ser. B}
\bvolume{47}
\bpages{238--252}.
\bid{issn={0035-9246}, mr={0816088}}
\end{barticle}
%
\endbibitem

\bibitem{gelman}
%
\begin{barticle}[auto:STB|2011/09/12|07:03:23]
\bauthor{\bsnm{Gelman},~\bfnm{A.}\binits{A.}},
\bauthor{\bsnm{Bois},~\bfnm{F.~Y.}\binits{F.Y.}} \AND
\bauthor{\bsnm{Jiang},~\bfnm{J.}\binits{J.}}
(\byear{1996}).
\btitle{Physiological pharmacokinetic analysis using population
modeling and informative prior distributions}.
\bjournal{J. Amer. Statist. Assoc.}
\bvolume{91}
\bpages{1400--1412}.
\end{barticle}
%
\endbibitem

\bibitem{girolami}
%
\begin{barticle}[mr]
\bauthor{\bsnm{Girolami},~\bfnm{Mark}\binits{M.}}
(\byear{2008}).
\btitle{Bayesian inference for differential equations}.
\bjournal{Theoret. Comput. Sci.}
\bvolume{408}
\bpages{4--16}.
\bid{doi={10.1016/j.tcs.2008.07.005}, issn={0304-3975}, mr={2460604}}
\end{barticle}
%
\endbibitem

\bibitem{goldstein}
%
\begin{barticle}[mr]
\bauthor{\bsnm{Goldstein},~\bfnm{Larry}\binits{L.}} \AND
\bauthor{\bsnm{Messer},~\bfnm{Karen}\binits{K.}}
(\byear{1992}).
\btitle{Optimal plug-in estimators for nonparametric functional estimation}.
\bjournal{Ann. Statist.}
\bvolume{20}
\bpages{1306--1328}.
\bid{doi={10.1214/aos/1176348770}, issn={0090-5364}, mr={1186251}}
\end{barticle}
%
\endbibitem

\bibitem{hairer}
%
\begin{bbook}[mr]
\bauthor{\bsnm{Hairer},~\bfnm{E.}\binits{E.}} \AND
\bauthor{\bsnm{Wanner},~\bfnm{G.}\binits{G.}}
(\byear{1996}).
\btitle{Solving Ordinary Differential Equations. {II}. Stiff and Differential-Algebraic Problems},
\bedition{2nd} ed.
\bseries{Springer Series in Computational Mathematics}
\bvolume{14}.
\baddress{Berlin}: \bpublisher{Springer}.
\bid{mr={1439506}}
\end{bbook}
%
\endbibitem

\bibitem{hall}
%
\begin{barticle}[mr]
\bauthor{\bsnm{Hall},~\bfnm{Peter}\binits{P.}} \AND
\bauthor{\bsnm{Marron},~\bfnm{J.~S.}\binits{J.S.}}
(\byear{1990}).
\btitle{On variance estimation in nonparametric regression}.
\bjournal{Biometrika}
\bvolume{77}
\bpages{415--419}.
\bid{doi={10.1093/biomet/77.2.415}, issn={0006-3444}, mr={1064818}}
\end{barticle}
%
\endbibitem

\bibitem{hemker}
%
\begin{bincollection}[auto:STB|2011/09/12|07:03:23]
\bauthor{\bsnm{Hemker},~\bfnm{P.~W.}\binits{P.W.}}
(\byear{1972}).
\btitle{Numerical methods for differential equations in system
simulation and
in parameter estimation}.
In \bbooktitle{Analysis and Simulation of Biochemical Systems}
(\beditor{\bfnm{H.~C.}\binits{H.C.}~\bsnm{Hemker}} \AND
\beditor{\bfnm{B.}\binits{B.}~\bsnm{Hess}}, eds.)
\bpages{59--80}.
 \baddress{Amsterdam}: \bpublisher{North Holland}.
\end{bincollection}
%
\endbibitem

\bibitem{hlavacek}
%
\begin{barticle}[pbm]
\bauthor{\bsnm{Hlavacek},~\bfnm{W.~S.}\binits{W.S.}} \AND
\bauthor{\bsnm{Savageau},~\bfnm{M.~A.}\binits{M.A.}}
(\byear{1996}).
\btitle{Rules for coupled expression of regulator and effector genes in
inducible circuits}.
\bjournal{J. Mol. Biol.}
\bvolume{255}
\bpages{121--139}.
\bid{doi={10.1006/jmbi.1996.0011}, issn={0022-2836},
pii={S0022-2836(96)90011-X}, pmid={8568860}}
\end{barticle}
%
\endbibitem

\bibitem{hoeffding}
%
\begin{barticle}[mr]
\bauthor{\bsnm{Hoeffding},~\bfnm{Wassily}\binits{W.}}
(\byear{1963}).
\btitle{Probability inequalities for sums of bounded random variables}.
\bjournal{J.~Amer. Statist. Assoc.}
\bvolume{58}
\bpages{13--30}.
\bid{issn={0162-1459}, mr={0144363}}
\end{barticle}
%
\endbibitem

\bibitem{hooker}
%
\begin{barticle}[mr]
\bauthor{\bsnm{Hooker},~\bfnm{Giles}\binits{G.}}
(\byear{2009}).
\btitle{Forcing function diagnostics for nonlinear dynamics}.
\bjournal{Biometrics}
\bvolume{65}
\bpages{928--936}.
\bid{doi={10.1111/j.1541-0420.2008.01172.x}, issn={0006-341X}, mr={2649866}}
\end{barticle}
%
\endbibitem

\bibitem{huber}
%
\begin{bbook}[mr]
\bauthor{\bsnm{Huber},~\bfnm{Peter~J.}\binits{P.J.}}
(\byear{1981}).
\btitle{Robust Statistics}.
\baddress{New York}: \bpublisher{Wiley}.
\bid{mr={0606374}}
\end{bbook}
%
\endbibitem

\bibitem{jennrich}
%
\begin{barticle}[mr]
\bauthor{\bsnm{Jennrich},~\bfnm{Robert~I.}\binits{R.I.}}
(\byear{1969}).
\btitle{Asymptotic properties of non-linear least squares estimators}.
\bjournal{Ann. Math. Statist.}
\bvolume{40}
\bpages{633--643}.
\bid{issn={0003-4851}, mr={0238419}}
\end{barticle}
%
\endbibitem

\bibitem{jones}
%
\begin{barticle}[mr]
\bauthor{\bsnm{Jones},~\bfnm{M.~C.}\binits{M.C.}},
\bauthor{\bsnm{Marron},~\bfnm{J.~S.}\binits{J.S.}} \AND
\bauthor{\bsnm{Sheather},~\bfnm{S.~J.}\binits{S.J.}}
(\byear{1996}).
\btitle{A brief survey of bandwidth selection for density estimation}.
\bjournal{J. Amer. Statist. Assoc.}
\bvolume{91}
\bpages{401--407}.
\bid{issn={0162-1459}, mr={1394097}}
\end{barticle}
%
\endbibitem

\bibitem{kikuchi}
%
\begin{barticle}[pbm]
\bauthor{\bsnm{Kikuchi},~\bfnm{Shinichi}\binits{S.}},
\bauthor{\bsnm{Tominaga},~\bfnm{Daisuke}\binits{D.}},
\bauthor{\bsnm{Arita},~\bfnm{Masanori}\binits{M.}},
\bauthor{\bsnm{Takahashi},~\bfnm{Katsutoshi}\binits{K.}} \AND
\bauthor{\bsnm{Tomita},~\bfnm{Masaru}\binits{M.}}
(\byear{2003}).
\btitle{Dynamic modeling of genetic networks using genetic algorithm and
S-system}.
\bjournal{Bioinformatics}
\bvolume{19}
\bpages{643--650}.
\bid{issn={1367-4803}, pmid={12651723}}
\end{barticle}
%
\endbibitem

\bibitem{liang}
%
\begin{barticle}[mr]
\bauthor{\bsnm{Liang},~\bfnm{Hua}\binits{H.}} \AND
\bauthor{\bsnm{Wu},~\bfnm{Hulin}\binits{H.}}
(\byear{2008}).
\btitle{Parameter estimation for differential equation models using a framework
of measurement error in regression models}.
\bjournal{J. Amer. Statist. Assoc.}
\bvolume{103}
\bpages{1570--1583}.
\bid{doi={10.1198/016214508000000797}, issn={0162-1459}, mr={2504205}}
\end{barticle}
%
\endbibitem

\bibitem{loader}
%
\begin{barticle}[mr]
\bauthor{\bsnm{Loader},~\bfnm{Clive~R.}\binits{C.R.}}
(\byear{1999}).
\btitle{Bandwidth selection: Classical or plug-in?}
\bjournal{Ann. Statist.}
\bvolume{27}
\bpages{415--438}.
\bid{doi={10.1214/aos/1018031201}, issn={0090-5364}, mr={1714723}}
\end{barticle}
%
\endbibitem

\bibitem{levenberg}
%
\begin{barticle}[mr]
\bauthor{\bsnm{Marquardt},~\bfnm{Donald~W.}\binits{D.W.}}
(\byear{1963}).
\btitle{An algorithm for least-squares estimation of nonlinear parameters}.
\bjournal{J. Soc. Indust. Appl. Math.}
\bvolume{11}
\bpages{431--441}.
\bid{mr={0153071}}
\end{barticle}
%
\endbibitem

\bibitem{mcmurry}
%
\begin{barticle}[mr]
\bauthor{\bsnm{McMurry},~\bfnm{Timothy~L.}\binits{T.L.}} \AND
\bauthor{\bsnm{Politis},~\bfnm{Dimitris~N.}\binits{D.N.}}
(\byear{2004}).
\btitle{Nonparametric regression with infinite order flat-top kernels}.
\bjournal{J. Nonparametr. Stat.}
\bvolume{16}
\bpages{549--562}.
\bid{doi={10.1080/10485250310001622596}, issn={1048-5252}, mr={2073041}}
\end{barticle}
%
\endbibitem

\bibitem{goldstein2}
%
\begin{barticle}[mr]
\bauthor{\bsnm{Messer},~\bfnm{Karen}\binits{K.}} \AND
\bauthor{\bsnm{Goldstein},~\bfnm{Larry}\binits{L.}}
(\byear{1993}).
\btitle{A new class of kernels for nonparametric curve estimation}.
\bjournal{Ann. Statist.}
\bvolume{21}
\bpages{179--195}.
\bid{doi={10.1214/aos/1176349021}, issn={0090-5364}, mr={1212172}}
\end{barticle}
%
\endbibitem

\bibitem{pollard}
%
\begin{barticle}[mr]
\bauthor{\bsnm{Pollard},~\bfnm{David}\binits{D.}} \AND
\bauthor{\bsnm{Radchenko},~\bfnm{Peter}\binits{P.}}
(\byear{2006}).
\btitle{Nonlinear least-squares estimation}.
\bjournal{J. Multivariate Anal.}
\bvolume{97}
\bpages{548--562}.
\bid{doi={10.1016/j.jmva.2005.04.002}, issn={0047-259X}, mr={2234037}}
\end{barticle}
%
\endbibitem

\bibitem{priestley}
%
\begin{barticle}[mr]
\bauthor{\bsnm{Priestley},~\bfnm{M.~B.}\binits{M.B.}} \AND
\bauthor{\bsnm{Chao},~\bfnm{M.~T.}\binits{M.T.}}
(\byear{1972}).
\btitle{Non-parametric function fitting}.
\bjournal{J. Roy. Statist. Soc. Ser. B}
\bvolume{34}
\bpages{385--392}.
\bid{issn={0035-9246}, mr={0331616}}
\end{barticle}
%
\endbibitem

\bibitem{qi}
%
\begin{barticle}[mr]
\bauthor{\bsnm{Qi},~\bfnm{Xin}\binits{X.}} \AND
\bauthor{\bsnm{Zhao},~\bfnm{Hongyu}\binits{H.}}
(\byear{2010}).
\btitle{Asymptotic efficiency and finite-sample properties of the generalized
profiling estimation of parameters in ordinary differential equations}.
\bjournal{Ann. Statist.}
\bvolume{38}
\bpages{435--481}.
\bid{doi={10.1214/09-AOS724}, issn={0090-5364}, mr={2589327}}
\end{barticle}
%
\endbibitem

\bibitem{ramsay}
%
\begin{barticle}[mr]
\bauthor{\bsnm{Ramsay},~\bfnm{J.~O.}\binits{J.O.}},
\bauthor{\bsnm{Hooker},~\bfnm{G.}\binits{G.}},
\bauthor{\bsnm{Campbell},~\bfnm{D.}\binits{D.}} \AND
\bauthor{\bsnm{Cao},~\bfnm{J.}\binits{J.}}
(\byear{2007}).
\btitle{Parameter estimation for differential equations: A generalized
smoothing approach}.
\bjournal{J. R. Stat. Soc. Ser. B Stat. Methodol.}
\bvolume{69}
\bpages{741--796}.
\bnote{With discussions and a reply by the authors}.
\bid{doi={10.1111/j.1467-9868.2007.00610.x}, issn={1369-7412}, mr={2368570}}
\end{barticle}
%
\endbibitem

\bibitem{schuster}
%
\begin{barticle}[mr]
\bauthor{\bsnm{Schuster},~\bfnm{E.}\binits{E.}} \AND
\bauthor{\bsnm{Yakowitz},~\bfnm{S.}\binits{S.}}
(\byear{1979}).
\btitle{Contributions to the theory of nonparametric regression, with
application to system identification}.
\bjournal{Ann. Statist.}
\bvolume{7}
\bpages{139--149}.
\bid{issn={0090-5364}, mr={0515689}}
\end{barticle}
%
\endbibitem

\bibitem{sontag}
%
\begin{barticle}[mr]
\bauthor{\bsnm{Sontag},~\bfnm{Eduardo~D.}\binits{E.D.}}
(\byear{2001}).
\btitle{Structure and stability of certain chemical networks and applications
to the kinetic proofreading model of {T}-cell receptor signal transduction}.
\bjournal{IEEE Trans. Automat. Control}
\bvolume{46}
\bpages{1028--1047}.
\bid{doi={10.1109/9.935056}, issn={0018-9286}, mr={1842137}}
\end{barticle}
%
\endbibitem

\bibitem{stigler}
%
\begin{barticle}[mr]
\bauthor{\bsnm{Stigler},~\bfnm{Stephen~M.}\binits{S.M.}}
(\byear{1981}).
\btitle{Gauss and the invention of least squares}.
\bjournal{Ann. Statist.}
\bvolume{9}
\bpages{465--474}.
\bid{issn={0090-5364}, mr={0615423}}
\end{barticle}
%
\endbibitem

\bibitem{stortelder}
%
\begin{barticle}[auto:STB|2011/09/12|07:03:23]
\bauthor{\bsnm{Stortelder},~\bfnm{W.~J.~H.}\binits{W.J.H.}}
(\byear{1996}).
\btitle{Parameter estimation in dynamic systems}.
\bjournal{Math. Comput. Simulat.}
\bvolume{42}
\bpages{135--142}.
\end{barticle}
%
\endbibitem

\bibitem{szego}
%
\begin{bbook}[mr]
\bauthor{\bsnm{Szeg{\H{o}}},~\bfnm{G{\'a}bor}\binits{G.}}
(\byear{1975}).
\btitle{Orthogonal Polynomials}, \bedition{4th} ed.
\bseries{American Mathematical Society, Colloquium Publications} \bvolume{XXIII}.
\baddress{Providence, RI}: \bpublisher{Amer. Math. Soc.}
\bid{mr={0372517}}
\end{bbook}
%
\endbibitem

\bibitem{tsybakov}
%
\begin{bbook}[mr]
\bauthor{\bsnm{Tsybakov},~\bfnm{Alexandre~B.}\binits{A.B.}}
(\byear{2009}).
\btitle{Introduction to Nonparametric Estimation}.
\bseries{Springer Series in Statistics}.
\baddress{New York}: \bpublisher{Springer}.
\bid{mr={2724359}}
\end{bbook}
%
\endbibitem

\bibitem{geer1}
%
\begin{barticle}[mr]
\bauthor{\bparticle{van~de} \bsnm{Geer},~\bfnm{Sara}\binits{S.}}
(\byear{1990}).
\btitle{Estimating a regression function}.
\bjournal{Ann. Statist.}
\bvolume{18}
\bpages{907--924}.
\bid{doi={10.1214/aos/1176347632}, issn={0090-5364}, mr={1056343}}
\end{barticle}
%
\endbibitem

\bibitem{geer2}
%
\begin{barticle}[mr]
\bauthor{\bparticle{van~de} \bsnm{Geer},~\bfnm{Sara}\binits{S.}}
\AND
\bauthor{\bsnm{Wegkamp},~\bfnm{Marten}\binits{M.}}
(\byear{1996}).
\btitle{Consistency for the least squares estimator in nonparametric
regression}.
\bjournal{Ann. Statist.}
\bvolume{24}
\bpages{2513--2523}.
\bid{doi={10.1214/aos/1032181165}, issn={0090-5364}, mr={1425964}}
\end{barticle}
%
\endbibitem

\bibitem{vdvaart}
%
\begin{bbook}[mr]
\bauthor{\bparticle{van~der} \bsnm{Vaart},~\bfnm{A.~W.}\binits{A.W.}}
(\byear{1998}).
\btitle{Asymptotic Statistics}.
\bseries{Cambridge Series in Statistical and Probabilistic Mathematics}
\bvolume{3}.
\baddress{Cambridge}: \bpublisher{Cambridge Univ. Press}.
\bid{mr={1652247}}
\end{bbook}
%
\endbibitem

\bibitem{wellner}
%
\begin{bbook}[mr]
\bauthor{\bparticle{van~der} \bsnm{Vaart},~\bfnm{Aad~W.}\binits
{A.W.}} \AND
\bauthor{\bsnm{Wellner},~\bfnm{Jon~A.}\binits{J.A.}}
(\byear{2000}).
\btitle{Weak Convergence and Empirical Processes: With Applications to Statistics},
\bedition{2}nd ed.
\bseries{Springer Series in Statistics}.
\baddress{New York}: \bpublisher{Springer}.
\end{bbook}
%
\endbibitem

\bibitem{vanes}
%
\begin{bbook}[mr]
\bauthor{\bparticle{van} \bsnm{Es},~\bfnm{A.~J.}\binits{A.J.}}
(\byear{1991}).
\btitle{Aspects of Nonparametric Density Estimation}.
\bseries{CWI Tract}
\bvolume{77}.
\baddress{Amsterdam}: \bpublisher{Stichting Mathematisch Centrum
Centrum voor
Wiskunde en Informatica}.
\bid{mr={1129890}}
\end{bbook}
%
\endbibitem

\bibitem{varah}
%
\begin{barticle}[mr]
\bauthor{\bsnm{Varah},~\bfnm{J.~M.}\binits{J.M.}}
(\byear{1982}).
\btitle{A spline least squares method for numerical parameter
estimation in
differential equations}.
\bjournal{SIAM J. Sci. Statist. Comput.}
\bvolume{3}
\bpages{28--46}.
\bid{doi={10.1137/0903003}, issn={0196-5204}, mr={0651865}}
\end{barticle}
%
\endbibitem

\bibitem{voit}
%
\begin{bbook}[auto:STB|2011/09/12|07:03:23]
\bauthor{\bsnm{Voit},~\bfnm{E.~O.}\binits{E.O.}}
(\byear{2000}).
\btitle{Computational Analysis of Biochemical Systems: A Practical
Guide for
Biochemists and Molecular Biologists}.
\baddress{Cambridge}: \bpublisher{Cambridge Univ. Press}.
\end{bbook}
%
\endbibitem

\bibitem{almeida}
%
\begin{barticle}[auto:STB|2011/09/12|07:03:23]
\bauthor{\bsnm{Voit},~\bfnm{E.~O.}\binits{E.O.}} \AND
\bauthor{\bsnm{Almeida},~\bfnm{J.}\binits{J.}}
(\byear{2004}).
\btitle{Decoupling dynamical systems for pathway identification from metabolic
profiles}.
\bjournal{Bioinformatics}
\bvolume{10}
\bpages{1670--1681}.
\end{barticle}
%
\endbibitem

\bibitem{savageau}
%
\begin{barticle}[auto:STB|2011/09/12|07:03:23]
\bauthor{\bsnm{Voit},~\bfnm{E.~O.}\binits{E.O.}} \AND
\bauthor{\bsnm{Savageau},~\bfnm{M.~A.}\binits{M.A.}}
(\byear{1982}).
\btitle{Power-law approach to modeling biological systems; III.
Methods of
analysis}.
\bjournal{J. Ferment. Technol.}
\bvolume{60}
\bpages{233--241}.
\end{barticle}
%
\endbibitem

\bibitem{biegler}
%
\begin{barticle}[mr]
\bauthor{\bsnm{W{\"a}chter},~\bfnm{Andreas}\binits{A.}} \AND
\bauthor{\bsnm{Biegler},~\bfnm{Lorenz~T.}\binits{L.T.}}
(\byear{2006}).
\btitle{On the implementation of an interior-point filter line-search algorithm
for large-scale nonlinear programming}.
\bjournal{Math. Program.}
\bvolume{106}
\bpages{25--57}.
\bid{doi={10.1007/s10107-004-0559-y}, issn={0025-5610}, mr={2195616}}
\end{barticle}
%
\endbibitem

\bibitem{wand2}
%
\begin{bbook}[mr]
\bauthor{\bsnm{Wand},~\bfnm{M.~P.}\binits{M.P.}} \AND
\bauthor{\bsnm{Jones},~\bfnm{M.~C.}\binits{M.C.}}
(\byear{1995}).
\btitle{Kernel Smoothing}.
\bseries{Monographs on Statistics and Applied Probability}
\bvolume{60}.
\baddress{London}: \bpublisher{Chapman \& Hall}.
\bid{mr={1319818}}
\end{bbook}
%
\endbibitem

\bibitem{schucany}
%
\begin{barticle}[mr]
\bauthor{\bsnm{Wand},~\bfnm{Matthew~P.}\binits{M.P.}} \AND
\bauthor{\bsnm{Schucany},~\bfnm{William~R.}\binits{W.R.}}
(\byear{1990}).
\btitle{Gaussian-based kernels}.
\bjournal{Canad. J. Statist.}
\bvolume{18}
\bpages{197--204}.
\bid{doi={10.2307/3315450}, issn={0319-5724}, mr={1079592}}
\end{barticle}
%
\endbibitem

\bibitem{wasserman}
%
\begin{bbook}[mr]
\bauthor{\bsnm{Wasserman},~\bfnm{Larry}\binits{L.}}
(\byear{2006}).
\btitle{All of Nonparametric Statistics}.
\bseries{Springer Texts in Statistics}.
\baddress{New York}: \bpublisher{Springer}.
\bid{mr={2172729}}
\end{bbook}
%
\endbibitem

\bibitem{mathematica}
%
\begin{bmisc}[auto:STB|2011/09/12|07:03:23]
\borganization{Wolfram Research, Inc.}
(\byear{2007}).
\bhowpublished{Mathematica, Version 6.0. Champaign, IL}.
\end{bmisc}
%
\endbibitem

\bibitem{wu}
%
\begin{barticle}[mr]
\bauthor{\bsnm{Wu},~\bfnm{Chien-Fu}\binits{C.F.}}
(\byear{1981}).
\btitle{Asymptotic theory of nonlinear least squares estimation}.
\bjournal{Ann. Statist.}
\bvolume{9}
\bpages{501--513}.
\bid{issn={0090-5364}, mr={0615427}}
\end{barticle}
%
\endbibitem

\bibitem{xue}
%
\begin{barticle}[mr]
\bauthor{\bsnm{Xue},~\bfnm{Hongqi}\binits{H.}},
\bauthor{\bsnm{Miao},~\bfnm{Hongyu}\binits{H.}} \AND
\bauthor{\bsnm{Wu},~\bfnm{Hulin}\binits{H.}}
(\byear{2010}).
\btitle{Sieve estimation of constant and time-varying coefficients in nonlinear
ordinary differential equation models by considering both numerical
error and
measurement error}.
\bjournal{Ann. Statist.}
\bvolume{38}
\bpages{2351--2387}.
\bid{doi={10.1214/09-AOS784}, issn={0090-5364}, mr={2676892}}
\end{barticle}
%
\endbibitem

\end{thebibliography}
\end{document}